\documentclass{amsart}
\usepackage[dvips]{graphicx}

\def\nt{\noindent}

\newtheorem{thm}{Theorem}[section]

\newtheorem{lem}[thm]{Lemma}
\newtheorem{prop}[thm]{Proposition}
\newtheorem{defn}[thm]{Definition}
\newtheorem{rem}[thm]{Remark}
\numberwithin{equation}{subsection}

\def\R{\mathbb{R}}
\def\S{\mathcal{S}}
\def\V{\mathbb{V}}
\def\C{\mathbb{C}}
\def\N{\mathbb{N}}
\def\reg{\mathrm{reg}}

\def\clas{\mathrm{cl}}
\def\sing{\mathrm{sing}}
\def\supp{\mathrm{supp}\hspace{.05cm}}
\def\Diff{\mathrm{Diff}}
\def\Op{\mathrm{Op}}
\def\op{\mathrm{op}}
\def\M{\mathbb M}
\def\O{\mathcal{O}}
\def\P{\mathcal P}
\def\Q{\mathcal Q}
\def\X{X^\wedge}
\def\nt{\noindent}
\def\Int{\displaystyle\int}
\def\TenH{\otimes_H}

\def\TenP{\widehat{\otimes}_\pi}
\def\TenG{\widehat{\otimes}_\Gamma}
\def\Pr{\noindent\textbf{Proof. }}
\def\eop{\hfill\fbox\\}
\def\Fr{Fr\'{e}chet }

\def\Kl{\kappa_\lambda}
\def\Kll{\set{\kappa_\lambda}_{\lambda\in\R_+}}
\def\Re{\mathrm{Re}\,}
\def\Im{\mathrm{Im}\,}
\def\RR{\mathcal{R}}
\def\WW{\mathcal W}
\def\MM{\mathcal M}
\def\VV{\mathcal{V}}
\def\U{\mathcal{U}}
\def\W{\mathbb{W}}
\def\C{\mathbb{C}}
\def\Z{\mathbb{Z}}

\def\A{\mathcal{A}}
\def\H{\mathcal{H}}
\def\L{\mathcal{L}}
\def\K{\mathcal{K}}
\def\Ksg{\mathcal{K}^{s,\gamma}}
\def\KsgX{\mathcal{K}^{s,\gamma}(X^\wedge)}
\def\E{\mathcal{E}}

\def\Summa{\displaystyle\sum}
\def\dim{\mathrm{dim}\,}
\def\Diff{\mathrm{Diff}}
\def\id{\mathrm{id}}
\def\clas{\mathrm{cl}}
\def\cone{\mathrm{cone}}
\def\loc{\mathrm{loc}}
\def\op{\mathrm{op}}
\def\Op{\mathrm{Op}}

\newcommand{\norm}[1]{\left\Vert#1\right\Vert}
\newcommand{\abs}[1]{\left\vert#1\right\vert}
\newcommand{\dslash}{d\llap {\raisebox{.9ex}{$\scriptstyle-\!$}}}
\newcommand{\Set}[1]{\left(#1\right)}
\newcommand{\dint}{\int\hspace{-.15cm}\int}
\newcommand{\set}[1]{\left\{#1\right\}}
\newcommand{\virg}[1]{`{#1}'}
\newcommand{\Rr}{\overline{\mathbb R}_+}
\newcommand{\Dr}{r\frac{\partial}{\partial r}}
\newcommand{\eps}{\varepsilon}
\newcommand{\To}{\longrightarrow}
\newcommand{\ang}[1]{\langle{#1}\rangle}
\newcommand{\til}[1]{\widetilde{#1}}
\newcommand{\cl}[1]{\overline{#1}}
\newcommand{\Projlim}[1]{\displaystyle\lim_{\substack{\longleftarrow\\ {#1}}}}
\newcommand{\KsgXAs}[1]{\mathcal{K}^{s,\gamma}_{#1}(X^\wedge)}

\begin{document}

\title{Green Operators in the Edge Calculus}%
\author{B.-W. Schulze and A.Volpato}%

\begin{abstract}
The task to construct parametrices of elliptic differential
operators on a manifold with edges requires a calculus of operators
with a two-component principal symbolic hierarchy, consisting of
(edge-degenerate) interior and (operator-valued) edge symbols. This
so-called edge-algebra can be interpreted as a generalisation of the
(pseudo-differential) algebra of boundary value problems without the
transmission property at the boundary. We study new properties of
the edge-algebra, in particular, what concerns the role of Green
operators and their kernel representations.
\end{abstract}
\maketitle

\nt\textbf{2000 AMS-classification:} 35S15, 47G30, 58J05\\

\nt\textbf{Keywords:} edge-degenerate operators in weighted spaces,
operators on manifolds with edges\\

\tableofcontents


\section*{Introduction}
\nt This paper is aimed at characterising the structure of Green
operators that occur in Green's functions of elliptic boundary value
problems (without or with the transmission property at the boundary)
or, more generally, in elliptic (so-called) edge problems. A
boundary or an edge represents a singularity of the configuration,
and this causes a typical singular behaviour of solutions to such
problems, here expressed in terms of discrete or continuous
asymptotics. The Green operators just encode such asymptotics. At
the same time it is interesting to understand the kernel structure
of trace and potential operators in the boundary (edge) calculus;
this is of a similar kind as the one of Green operators.\\

\nt By a manifold $M$ with edge $Y$ we understand a topological
space such that $M\setminus Y$ and $Y$ are $C^\infty$ manifolds, and
$Y$ has a neighbourhood $V$ in $M$ equipped with the structure of an
$X^\triangle$-bundle over $Y$ for a closed compact $C^\infty$
manifold $X;$ here $X^\triangle:=(\Rr\times X)/(\{0\}\times X)$ is
the infinite cone with base $X.$ In addition $V$ is the quotient
space of an $\Rr\times X$-bundle $\V$ over $Y$ where $\V\rightarrow
V$ is defined by the fibrewise projection $\Rr\times X\rightarrow
X^\triangle.$ From the definition it follows that $M$ itself is the
quotient space of a $C^\infty$ manifold $\M$ with boundary
$\partial\M,$ and $\partial\M$ is an $X$-bundle over $Y.$ Locally
near $Y$ the manifold $M$ has the structure of a wedge
$X^\triangle\times\Omega,$ $\Omega\subseteq\R^q$ open, $q=\dim Y,$
and $\M,$ the so-called stretched manifold of $M,$ is locally near
$\partial\M$ of the form $\Rr\times X\times\Omega.$ The starting
point of our discussion are differential operators on
$\M_\reg:=\M\setminus\partial\M$ with smooth coefficients that are
edge-degenerate near $\partial\M,$ i.e., locally, in the splitting
of variables $(r,x,y)\in\R_+\times X\times\Omega$ of the form
\begin{equation}\label{neu}
A=r^{-\mu}\sum_{j+\abs\alpha\leq\mu}a_{j\alpha}(r,y)\Set{-\Dr}^j\Set{rD_y}^\alpha
\end{equation}
with coefficients $a_{j\alpha}\in
C^\infty(\Rr\times\Omega,\Diff^{\mu-(j+\abs\alpha)}(X)).$ Here
$\Diff^\nu(\cdot)$ denotes the space of all differential operators
of order $\nu$ on the manifold in parentheses. If $A$ is elliptic in
the sense of standard ellipticity in $\M_\reg,$ together with the
condition that for every fixed $y$ the expression
$$\sum_{j+\abs\alpha\leq\mu}a_{j\alpha}(0,y)(-i\rho)^j\eta^\alpha$$
represents a parameter-dependent elliptic family of differential
operators on $X$ with the paramenters $(\rho,\eta)\in\R^{1+q}$ (in
the sense of Agranovich and Vishik, cf. \cite{Agra1}), we may ask
the (pseudo-differential) nature of a parametrix of $A$ and the
Fredholm property in suitable weighted Sobolev spaces on $\M.$ It
turns out that such a characterisation of the solvability properties
of the equation $Au=f$ induces a large variety of additional
structures of an \virg{unexpected} complexity.\\

\nt Comparing the situation with boundary value problems (taking
into account that a, say, compact, $C^\infty$ manifold $M$ with
boundary is a particular manifold with edge, where the above
mentioned cone bundle near the edge $\partial M$ can be identified
with the inner normal bundle) we see that the solvability should be
discussed in connection with extra conditions of trace (and, in
general, also potential) type with respect to the edge, analogously
as elliptic boundary conditions. In addition the edge-degenerate
nature of the operators makes it necessary to refer to specific
weighted edge spaces rather than \virg{standard} Sobolev spaces.
This aspect requires some special attention, and there are many
possible choices of such spaces. Answers of the above mentioned
solvability problem (together with the construction of the
additional edge conditions) for a certain special class of such
spaces, including the characterisation of the structure of
parametrices, may be found in \cite{Schu32}, \cite{Schu2}.\\

\nt Analogously as for the case with boundary, in the case of a
manifold with edge we talk about an edge problem for $A.$ To be more
precise, by an edge problem we understand the construction (or
characterisation) of a solution $u$ to the equation $Au=f$ on a
manifold with edge $Y$ under suitable trace (or edge) conditions on
$u$ at $Y$ and, if necessary, additional potentials $Kv$ of
distributions $v$ on $Y$ which are to be added to $Au$ in the image
of $A.$\\

\nt Formally, Green operators have an analogous origin as Green's
functions in the parametrices of elliptic boundary value problems
(Green's function has the form of a sum $E+G$ where $E$ is a
fundamental solution or a parametrix of the given elliptic operator
and $G$ a Green operator in our sense). However, in the edge case
(or when the transmission property is violated) the Green operators
are \virg{loaded} with a huge variety of asymptotic information at
the edge, coming from the zeros of the conormal symbol of the given
operator, see the formula (\ref{cosy}) below. In the special case of
boundary value problems with the transmission property the
asymptotics are nothing other than smoothness up to the boundary,
see, for instance, \cite{Bout2}, and kernel characterisations of
corresponding Green symbols may be found in \cite{Schu16}. Kernel
characterisations for Green operators of the cone calculus are given
in \cite{Seil2}; the case of Green edge symbols is treated in
\cite{Schu13}, although this refers to projective tensor products of
spaces with discrete asymptotics; here we take another (completely
adequate but more convenient) tensor product and also study the case
of continuous asymptotics, see also \cite{Volp}. The discrete case
makes sense for $y$-independent asymptotics. In the $y$-dependent
case it is reasonable to formulate the phenomena in the frame of
continuous asymptotics (details will be explained below). Green
operators are smoothing on $\M_\reg.$ Close to
$\M_\sing:=\partial\M$ they can be characterised (modulo a certain
kind of global smoothing operators of the edge calculus) as
pseudo-differential operators on the edge with (classical)
operator-valued symbols, with \virg{twisted} homogeneity. Green
operators formally appear in a pseudo-differential calculus
containing the operator $A$ together with its parametrix in the
elliptic case. To illustrate that a little more we write
$$A=\Op(a)=F^{-1}_{\eta\rightarrow y}a(y,\eta)F_{y\rightarrow\eta},$$
(also denoted by $\Op_y(a)$) with the Fourier transform $F$ in
$y\in\R^q,$ and
\begin{equation}\label{a}
a(y,\eta):=r^{-\mu}\sum_{j+\abs\alpha\leq\mu}a_{j\alpha}(r,y)\Set{-\Dr}^j(r\eta)^\alpha.
\end{equation}

\nt Here, $a(y,\eta)$ is interpreted as an operator-valued amplitude
function, operating in suitable weighted distribution spaces on the
infinite (open stretched) cone $\X:=\R_+\times X$ with base $X.$
Since the specific information is coming from a neighbourhood of
$r=0,$ we may (and will) assume that the coefficients
$a_{j\alpha}(r,y)$ are independent of $r$ for $r>R$ for some $R>0.$
The behaviour of $A$ near the edge $Y$ very much depends on the
so-called principal edge symbol
$$\sigma_\wedge(A)(y,\eta):=r^{-\mu}\sum_{j+\abs\alpha\leq\mu}a_{j\alpha}(0,y)\Set{-\Dr}^j(r\eta)^\alpha,$$
$(y,\eta)\in T^*Y\setminus0,$ and the subordinate conormal symbol
\begin{equation}\label{cosy}
\sigma_c\sigma_\wedge(A)(y,z):=\sum^\mu_{j=0}a_{j0}(0,y)z^j
\end{equation}
which is (in the present case of a differential operator) an entire
function in $z\in\C$ taking values in $\Diff^\mu(X).$ As such it is
a holomorphic family of Fredholm operators
\begin{equation}\label{sigmac}
\sigma_c\sigma_\wedge(A)(y,z):H^s(X)\To H^{s-\mu}(X)
\end{equation}
in standard Sobolev spaces on $X.$ Parametrices of $A$ are to a
large extent determined by the inverse of $(\ref{sigmac})$ which is
a meromorphic operator function. Poles and multiplicities may depend
on $y;$ this gives rise to clouds of points in the complex plane,
and the $y$-dependent nature of meromorphy can be described in terms
of (operator-valued) analytic functionals in the complex plane,
pointwise discrete and branching, but smoothly depending on the
edge-variable $y$ (with respect to the Fr\'{e}chet topology of the
space of analytic functionals). Now the Green operators of the edge
calculus which participate in the parametrices of elliptic operators
inherit a part of this structure and affect the elliptic regularity
of solutions. To briefly recall the idea, an equation $Au=f$ with
elliptic $A$ is multiplied from the left by a parametrix $P$ of $A;$
then using $PA=I-G$ for a Green operator $G$ we obtain
\begin{equation}\label{Param}
PAu=f\quad\Rightarrow\quad u=Pf+Gu.
\end{equation}

\nt If $P$ is sensitive enough (which is the case in the
corresponding edge calculus) to transform a given $f$ with
asymptotics to a function $Pf$ with some resulting asymptotics and
if $G$ produces from any weighted distribution $u$ (without
asymptotics) a smooth function with asymptotics, from (\ref{Param})
we immediately conclude the asymptotics of $u.$ Note that a similar
philosophy works in every calculus with parametrices of elliptic
elements (with or without asymptotics), for instance, in the edge
calculus of \cite{Cori2}, where the base spaces of model cones may
have non-constant dimension.\\

\nt In other words, it is interesting to understand the nature of
Green operators in the edge calculus. Let us also note that the
above-mentioned trace and potential operators with respect to the
edge can also be subsumed under the concept of Green operators.


\section{Edge symbols and weighted Sobolev spaces}
{\setlength{\baselineskip}%
{.8\baselineskip} {\footnotesize \nt We study Mellin operators in
the distance variable $r\in\R_+$ to the singularity and recall the
edge quantisation of edge-degenerate symbols in $r$-direction. After
that we formulate the concept of discrete and continuous asymptotics
for $r\rightarrow0$ in weighted cone Sobolev spaces. Green symbols
of the edge calculus are formulated as classical symbols with
twisted homogeneity, mapping weighted distributions into smooth
functions with such asymptotics. } \par}


\subsection{Operator functions on the infinite stretched cone}\label{1.01}
Let $X$ be a closed compact $C^\infty$ manifold, and let
$L^\mu_\clas(X;\R^l)$ denote the space of classical
parameter-dependent pseudo-differential operators on $X$ of order
$\mu\in\R,$ with the parameter $\lambda\in\R^l$ (that is, the local
amplitude functions are classical symbols in
$(\xi,\lambda)\in\R^{n+l},$ $n=\dim X,$ and
$L^{-\infty}(X;\R^l)=\S(\R^l,L^{-\infty}(X)),$ where
$L^{-\infty}(X)$ is the space of smoothing operators on $X,$
identified with the space of operators with kernels in
$C^\infty(X\times X)$ via a Riemannian metric on $X).$ Analogously
as the operator functions $(\ref{a})$ in the pseudo-differential
case we consider
$$r^{-\mu}\op_r(p)(y,\eta)$$
for an operator family $p$ of the form
$$p(r,y,\rho,\eta):=\tilde p(r,y,r\rho,r\eta)$$
where $\tilde p(r,y,\tilde\rho,\tilde\eta)\in
C^\infty(\Rr\times\Omega,L^\mu_\clas(X;\R^{1+q}_{\tilde\rho,\tilde\eta}))$
and
$$\op_r(p)(y,\eta)u(r)=\dint e^{i(r-r')\rho}p(r,y,\rho,\eta)u(r')dr'\dslash\rho.$$

\nt Let $M$ denote the Mellin transform on $\R_+,$ i.e.,
$Mu(z):=\int_0^\infty r^{z-1}u(r)dr,$ first for $u\in
C_0^\infty(\R_+)$ and then extended to more general function and
distribution spaces, also vector-valued ones. Concerning details on
the Mellin transform in connection with edge-degenerate operators,
see \cite{Schu2} or \cite{Dors9}. Set $\Gamma_\delta:=\{z\in\C:\Re
z=\delta\}$ for any $\delta\in\R,$ and define the weighted Mellin
transform with weight $\gamma$ as the map $M_\gamma:u\mapsto
Mu(z)|_{\Gamma_{\frac12-\gamma}},$ $u\in C^\infty_0(\R_+).$ With
$M_\gamma$ we can associate weighted Mellin pseudo-differential
operators $\op_M^\gamma(f)$ with scalar as well as operator-valued
amplitude functions $f.$ Note that Mellin operators have been widely
used in Eskin's monograph \cite{Eski2} in connection with boundary
value problems without the transmission property, see also
\cite{Remp1}. In the edge calculus we take parameter-dependent
families of pseudo-differential operators $f(r,y,z,\eta)\in
C^\infty(\R_+\times\Omega,L^\mu_\clas(X;\Gamma_{\frac12-\gamma}\times\R^q));$
here $(z,\eta)\in\Gamma_{\frac12-\gamma}\times\R^q$ plays the role
of parameters while $(r,y)$ are additional variables. The associated
operators have the form
\begin{equation}\label{MellinOp}
\op^\gamma_M(f)(y,\eta)u(r):=\dint\Set{\frac{r}{r'}}^{-(\frac12-\gamma+i\rho)}
f(r,y,\frac12-\gamma+i\rho,\eta)u(r')\frac{dr'}{r'}\dslash\rho.
\end{equation}

\nt In the example of $(\ref{a})$ we can write
$$\tilde p(r,y,\tilde\rho,\tilde\eta)=\sum_{j+\abs\alpha\leq\mu}a_{j\alpha}(r,y)
(-i\tilde\rho)^j\tilde\eta^\alpha;$$ then
$a(y,\eta)=r^{-\mu}\op_r(p)(y,\eta),$ and
$$a(y,\eta):=r^{-\mu}\op^\beta_M(f)(y,\eta)$$
for $f(r,y,z,\eta):=\tilde f(r,y,z,r\eta),$ $\tilde
f(r,y,z,\tilde\eta)=\sum_{j+\abs\alpha\leq\mu}a_{j\alpha}(r,y)z^j\tilde\eta^\alpha$
and any $\beta\in\R$ (at this moment $a(y,\eta)$ is interpreted as
an operator function in the sense $C^\infty_0(\X)\To C^\infty(\X);$
later on we deal with weighted spaces on $\X=\R_+\times X).$\\

\nt Observe that the correspondence $p\rightarrow f$ in the case of
a differential operator (\ref{neu}) is canonical, i.e., we have in
this case
$$A=r^{-\mu}\Op_y(\op_r(p))=r^{-\mu}\Op_y(\op_M^\gamma(f)),$$
for any real $\gamma,$ or, alternatively,
\begin{equation}\label{MellQuant}
\op_r(p)(y,\eta)=\op_M^\gamma(f)(y,\eta)
\end{equation}
as elements of $C^\infty(\Omega,L^\mu_\clas(X^\wedge;\R^q)).$
However, in the pseudo-differential case the correspondence
$p\rightarrow f$ has to be achieved by a more subtle quantisation
which leaves smoothing remainders that are not necessarily
vanishing. It is known that such quantisations exist, indeed, see
\cite{Schu2}, and that for every $p$ there are elements $\tilde
f(r,y,z,\tilde\eta)\in
C^\infty(\Rr\times\Omega,L^\mu_\clas(X;\Gamma_{\frac12-\gamma}\times\R^q))$
such that for the associated $f$ the relation (\ref{MellQuant})
holds modulo such remainders. We even find holomorphic functions in
$z,$ see also Section \ref{3.2}, and (\ref{MellQuant}) then holds
for arbitrary $\gamma\in\R.$ The pseudo-differential edge calculus
of \cite{Schu2} contains amplitude functions of the form
$$\begin{array}{ll}
a(y,\eta)\hspace{-.3cm}&:=r^{-\mu}\sigma\{\omega(r[\eta])\op_M^{\gamma-\frac{n}2}(f)(y,\eta)\til\omega(r[\eta])\\
\\
&+(1-\omega(r[\eta]))\op_r(p)(y,\eta)(1-\til{\til\omega}(r[\eta]))\}\til\sigma,
\end{array}$$

\nt $n=\dim X,$ where $\omega,\til\omega,\til{\til\omega}$ and
$\sigma,\til\sigma$ are cut-off functions on the half axis, i.e.,
elements in $C_0^\infty(\Rr)$ which are equal to one near zero, and
with the property $\til{\til\omega}\prec\omega\prec\til\omega,$
where $\varphi\prec\psi$ means that $\psi\equiv1$ in a neighbourhood
of $\supp\varphi.$ In addition $\eta\rightarrow[\eta]$ is a strictly
positive $C^\infty$ function in $\eta\in\R^q$ such that
$[\eta]=\abs\eta$ for $\abs\eta>C$ for some $C>0.$\\

\nt In order to complete operators $\Op_y(a)$ to an algebra we have
to add other operator-valued symbols, mainly (so-called) Green
symbols and also smoothing Mellin ope-rator-families, see Section
\ref{3.2} below. Those are encoding the asymptotic properties of
solutions to elliptic equations when they are coming from the
parametrix.


\subsection{Weighted spaces with discrete asymptotics on cones}\label{1.1}

Green operators in the context of classical boundary value problems
are (locally in a collar neighbourhood of the boundary)
pseudo-differential operators along the boundary with symbols acting
as operators normal to the boundary. More precisely, the values of
the symbols are operators $G$ in $L^2(\R_+)$ such that
$G,G^*:L^2(\R_+)\To\S(\Rr)(=\S(\R)|_{\Rr})$ are continuous (this
concerns the so called type zero; otherwise the operators are
combined with differentiations transversal to the boundary). In the
generalisation to the case of a manifold with edges we replace the
inner normal $\Rr$ by a non-trivial model cone
$X^\bigtriangleup:=(\Rr\times X)/(\{0\}\times X)$ belonging to
corresponding local wedges, and the spaces $L^2(\R_+)$ and $\S(\Rr)$
by weighted spaces $\K^{0,\gamma}(X^\wedge)$ and
$\S^\beta_\P(X^\wedge),$ respectively, on the open stretched cone
$X^\wedge:=\R_+\times X\ni(r,x)$ with a certain behaviour for
$r\rightarrow0,$ encoded by a so called asymptotic type $\P.$ In
simplest cases asymptotics will have the form
\begin{equation}\label{new}
u(r,x)\sim\Summa_j\Summa^{m_j}_{k=0}c_{jk}(x)r^{-p_j}\log^kr
\qquad\mathrm{for}\quad r\rightarrow0
\end{equation}
with a sequence of triples $\P:=\{(p_j,m_j,L_j)\}_{j=0,1,\dots,N},$
$N\in\N\cup\{\infty\},$ $p_j\in\C,$ $m_j\in\N,$ and
finite-dimensional subspaces $L_j\subset C^\infty(X),$ such that
$c_{jk}\in L_j$ for all $0\le k\le m_j,$ and all $j.$ For $\dim X=0$
we have, in particular, a natural identification of $\S(\Rr)$ with
$\S^0_T(\R_+)$ for the Taylor asymptotic type
$T=\{(-j,0)\}_{j=0,1,\dots}$ (the spaces $L_j$ disappear in this
case). Let us now pass to the
precise definitions.\\

\nt We say that $\P$ is associated with weight data
$(\gamma,\Theta)$ for a weight $\gamma\in\R$ and
$\Theta=(\vartheta,0]$ for some $-\infty\le\vartheta<0$ if the set
$\pi_\C\P=\set{p_j}_{0\le j\le N}$ is contained in the strip
$\set{\frac{n+1}{2}-\gamma+\vartheta<\Re
z<\frac{n+1}{2}-\gamma},n=\dim X,$ $\pi_\C\P$ finite for finite
$\vartheta,$ and $\Re p_j\rightarrow-\infty$ as $j\rightarrow\infty$
when $\vartheta=-\infty$ and $N=\infty.$\\

\nt Given $\P$ associated with $(\gamma,\Theta)$ for a finite weight
interval $\Theta$ we set
$${\E}_{\P}(X^\wedge):=\Big\{\Summa^N_{j=0}\Summa^{m_j}_{k=0}
\omega(r)c_{jk}r^{-p_j}\log^kr:\textrm{$c_{jk}\in L_j$ for $0\le
k\le m_j,$ $0\le j\le N$}\Big\}.$$

\nt Moreover, let $\H^{s,\gamma}(X^\wedge)$ for $s\in\N,$
$\gamma\in\R,$ denote the subspace of all $u(r,x)\in
r^{\gamma-\frac{n}{2}}L^2(X^\wedge)$ (with $L^2$ referring to
$drdx$) such that
$$\Set{r\partial_r}^kD^\alpha_xu(r,x)\in r^{\gamma-\frac{n}{2}}L^2(X^\wedge)$$
for every $k\in\N,$ $\alpha\in\N^n,$ $k+\abs\alpha\le s;$ here
$D^\alpha_x:=v_1^{\alpha_1}\cdot$
$\hspace{-.1cm}\dots\hspace{-.1cm}$
$\cdot\hspace{.09cm}v_n^{\alpha_n}$ means the differentiation with
arbitrary vector fields $v_j$ on $X.$ In particular, we have
$\H^{0,0}(X^\wedge)=r^{-\frac{n}{2}}L^2(X^\wedge).$ Then duality and
interpolation give us a definition of
$\H^{s,\gamma}(X^\wedge)$ for arbitrary $s\in\R.$\\

\nt There is another useful scale of weighted spaces on $X^\wedge$
defined by
$$\KsgX:=\set{\omega f+(1-\omega)g:
f\in\H^{s,\gamma}(X^\wedge),~g\in H^s_\cone(X^\wedge)}$$ for some
cut-off function $\omega;$ see, for instance, \cite{Remp3} or
\cite{Schu2}. The space $H^s_\cone(X^\wedge)$ is defined to be the
set of all $g\in H^s_\loc(\R\times X)|_{\R_+\times X}$ such that for
every chart $\chi:\U\rightarrow B$ on $X$ to
$B:=\set{x\in\R^n:\abs{x}<1}$ and every $\varphi\in C^\infty_0(\U)$
we have
$$(1-\omega)\varphi g\in\Set{\beta\circ(1\times\chi)}^*H^s(\R^{n+1})|_\Gamma$$
for $\Gamma:=\set{(r,rx)\in\R^{n+1}:r\in\R_+,x\in B},$ $n=\dim X,$
and $\beta:\R_+\times B\rightarrow\Gamma,$ $\beta(r,x):=(r,rx),$
$(1\times\chi)(r,\cdot):=(r,\chi(\cdot)),$ see picture below.

\begin{figure}[h]
  \includegraphics[scale=0.5]{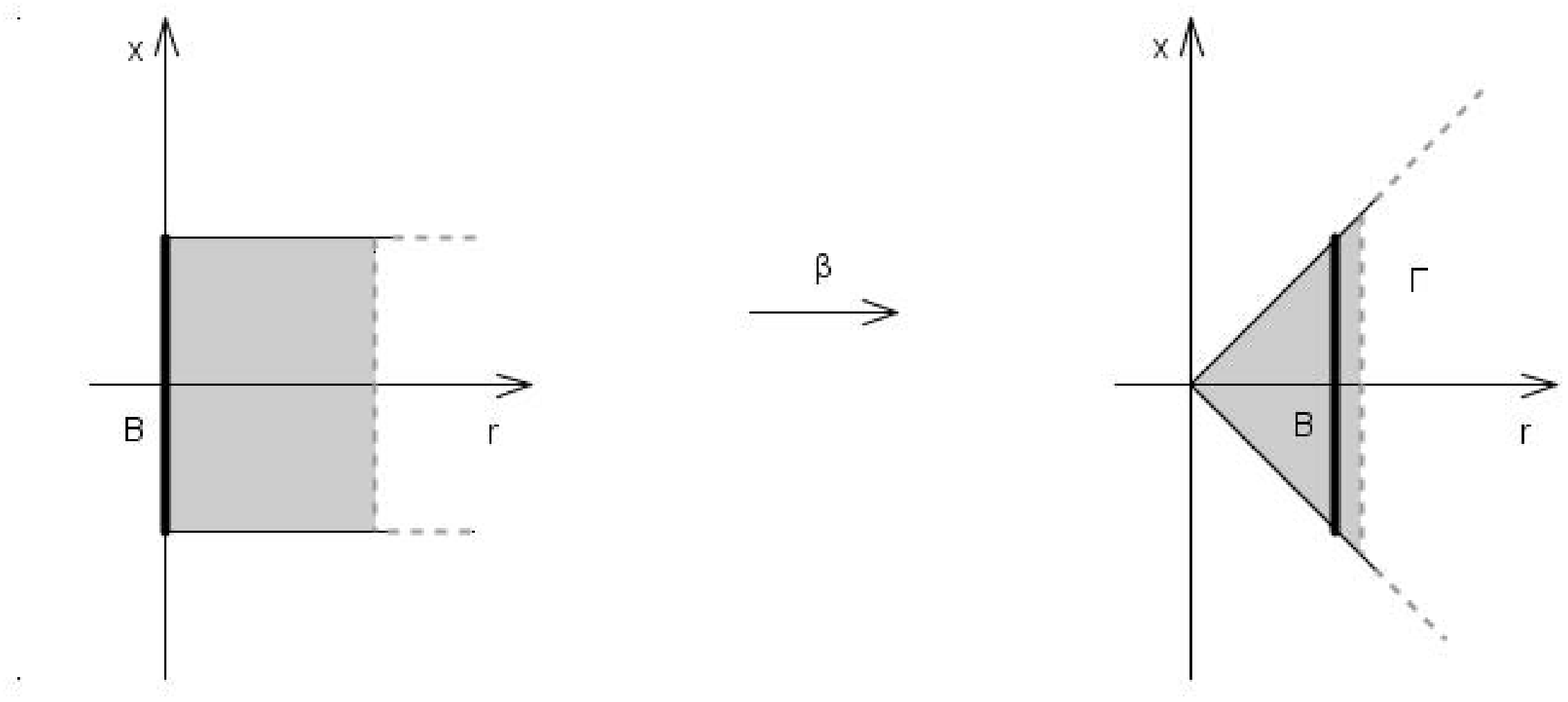}
\end{figure}

\nt The spaces $\KsgX$ can be endowed with scalar products such that
they are Hilbert spaces in a natural way; in particular,
$\K^{0,0}(X^\wedge)=\H^{0,0}(X^\wedge)=r^{-\frac{n}{2}}L^2(X^\wedge).$\\

\nt For a finite weight interval $\Theta=(\vartheta,0]$ we set
\begin{equation}\label{KTheta}
\Ksg_\Theta(X^\wedge):=\Projlim{N\in\N}\K^{s,\gamma-\vartheta-\frac{1}{N+1}}(X^\wedge)
\end{equation}
which is a \Fr space in the projective limit topology, and
$$\Ksg_\P(X^\wedge):=\Ksg_\Theta(X^\wedge)+\E_\P(X^\wedge),$$
as a direct sum, for every asymptotic type $\P$ which is associated
with the weight data $(\gamma,\Theta).$ For purposes below for every
$N\in\N$ and for $\gamma=0,$ we now form the spaces
\begin{equation}\label{BN}
B^N:=\ang{r}^{-N}\K^{N,0}(X^\wedge)
\end{equation}
and
\begin{equation}\label{APN}
A_\P^N:=\ang{r}^{-N}\K^{N,-\vartheta-\frac1{N+1}}(X^\wedge)+\E_\P(X^\wedge)
=w^{-\vartheta-\frac1{N+1}}B^N+\E_\P(X^\wedge);
\end{equation}
here, $w(r):=1+(r-1)\omega(r).$\\

\nt These are Hilbert spaces in a natural way, and we set
$$\S^0(X^\wedge):=\Projlim{N\in\N}B^N,\quad\S^0_\P(X^\wedge):=\Projlim{N\in\N}A^N_\P.$$

\nt More generally, we can form
$$\S^\gamma(X^\wedge):=w^\gamma\S^0(X^\wedge)
\quad\mathrm{and}\quad\S^\gamma_{T^{-\gamma}\P}(X^\wedge):=w^\gamma\S^0_\P(X^\wedge),$$

\nt where
$T^{-\gamma}\P:=\set{(p_j-\gamma,m_j,L_j)}_{j=0,\dots,N}.$\\

\begin{rem}\label{kappa}
\begin{itemize}
\item[$\mathrm{(i)}$] There are canonical continuous embeddings
\begin{equation}\label{Emb}
A^N_\P\hookrightarrow A^{N-1}_\P,\quad B^N\hookrightarrow B^{N-1}
\end{equation}
for all $N\ge1;$\\
\item[$\mathrm{(ii)}$] let us set
\begin{equation}\label{kl}
(\Kl u)(r,x):=\lambda^\frac{n+1}2u(\lambda r,x),
\end{equation}
$n=\dim X,$ $\lambda\in\R_+.$ Then we obtain strongly continuous
groups of isomorphisms
\begin{equation*}
\Kl:A^N_\P\rightarrow A^N_\P\quad as\;well\;as\quad
\Kl:B^N\rightarrow B^N
\end{equation*}
for every $N\in\N$ $($recall that $\Kll$ is said to be strongly
continuous on a Banach space $B$ if $\lambda\mapsto\Kl b$ represents
a continuous function $\R_+\rightarrow B$ for every $b\in B).$
\end{itemize}
\end{rem}


\subsection{Green edge symbols}\label{1.2}

Green symbols will be particular operator-valued symbols within the
framework of `twisted homogeneity'. Homogeneity in that sense means
the following. Let $E$ be a Hilbert space equipped with a strongly
continuous group $\{\Kl\}_{\lambda\in\R_+}$ of isomorphisms
$\Kl:E\rightarrow E,$ $\lambda\in\R_+,$ such that
$\Kl\kappa_\rho=\kappa_{\lambda\rho}$ for all $\lambda,\rho\in\R_+$
(in such a case we simply say that $E$ is endowed with a group
action). If $\til{E}$ is another Hilbert space with group action
$\{\til\kappa_\lambda\}_{\lambda\in\R_+},$ a $C^\infty$ function
$a_{(\mu)}(y,\eta)$ in $\Omega\times(\R^q\setminus\{0\}),$
$\Omega\subseteq\R^q$ open, with values in $\L(E,\til E)$ is called
homogeneous of order $\mu\in\R$ if
$$a_{(\mu)}(y,\lambda\eta)=\lambda^\mu\til\kappa_\lambda a_{(\mu)}(y,\eta)\Kl^{-1}$$
for all $\lambda\in\R_+.$\\

\nt Let us give the definition of the space of symbols
$S^\mu(\Omega\times\R^q;E,\til E).$ This space consists of the set
of all $C^\infty$ functions $a(y,\eta)$ in $\Omega\times\R^q$ with
values in $\L(E,\til E)$ such that
\begin{equation}\label{semi}
\sup_{\substack{y\in K \\
\eta\in\R^q}}\ang{\eta}^{|\beta|-\mu}\norm{\til\kappa_{\ang\eta}^{-1}\set{D^\alpha_yD^\beta_\eta
a(y,\eta)}\kappa_{\ang\eta}}_{\L(E,\til E)}
\end{equation}
is finite for every $K\subset\subset\Omega$ and every
$\alpha,\beta\in\N^q;$ here $\ang\eta:=(1+\abs\eta^2)^\frac12.$
Moreover, we denote by $S^{-\infty}(\Omega\times\R^q;E,\til
E):=\bigcap_{\mu\in\R} S^\mu(\Omega\times\R^q;E,\til E)$ the space
of symbols of order minus infinity.\\

\nt Symbols of that kind form a \Fr space with the expressions
(\ref{semi}) as semi-norms. They are `twisted' analogues of
H\"ormander's symbol spaces from the scalar case (i.e., when $E=\til
E=\C,$ $\Kl=\til\kappa_\lambda=\id_\C$ for all $\lambda\in\R_+$).
Standard manipulations known from the scalar case also make sense in
analogous form in the operator-valued case. In particular, we can
form asymptotic sums of sequences $a_j(y,\eta)$ of symbols the order
of which tend to $-\infty$ as $j\to\infty.$ Now
\begin{equation}\label{cl}
S^\mu_\clas(\Omega\times\R^q;E,\til E)
\end{equation}

\nt is defined as the subspace of $a(y,\eta)\in
S^\mu(\Omega\times\R^q;E,\til E)$ which admit asymptotic expansions
into symbols of the kind $\chi(\eta)a_{(\mu-j)}(y,\eta),$ $j\in\N,$
where $a_{(\mu-j)}(y,\eta)$ is homogeneous in the above sense, of
order $\mu-j$ and $\chi(\eta)$ an excision function, i.e., any
$\chi\in C^\infty(\R^q)$ that vanishes near $\eta=0$ and is equal to
$1$ for $\abs{\eta}\ge C$ for some $C>0.$
\begin{rem} Let
$\chi(\eta)$ be an excision function, and let $a_{(\mu)}(y,\eta)$ be
homogeneous of order $\mu$ as above; then
$$a(y,\eta):=\chi(\eta)a_{(\mu)}(y,\eta)\in S^\mu_\clas(\Omega\times\R^q;E,\til E).$$
\end{rem}

\nt In symbol spaces of the kind (\ref{cl}) for $E$ and $\til E$ we
will take, for instance,
\begin{equation}\label{Ksgamma}
\KsgX\quad\mathrm{and}\quad\S^\beta_\P(X^\wedge),
\end{equation}
respectively, for some discrete asymptotic type $\P$ (associated to
weight data $(\beta,\Theta),$ cf. Section \ref{1.1}). The spaces
(\ref{Ksgamma}) will be considered with the group action
(\ref{kl}).\\

\nt In order to define Green symbols we need a slight generalisation
to the case of \Fr spaces. We say, that a \Fr space $E,$ written as
the projective limit of a sequence of Hilbert spaces $E^j,$
$j\in\N,$ with continuous embedding $\dots\hookrightarrow
E^{j+1}\hookrightarrow E^j\hookrightarrow\dots\hookrightarrow E^0$
for all $j,$ is endowed with a group action
$\set{\Kl}_{\lambda\in\R_+},$ if $\set{\Kl}_{\lambda\in\R_+}$ is a
group action on $E^0$ and $\set{\Kl|_{E^j}}_{\lambda\in\R_+}$
defines a group action on $E^j$ for every $j.$

\begin{defn}\label{Green symbol}
An operator function $g(y,\eta)\in
C^\infty\big(\Omega\times\R^q,\L\big(\K^{0,\gamma}(X^\wedge),\linebreak
\K^{0,\beta}(X^\wedge)\big)\big)$ is said to be a Green symbol of
order $\mu\in\R,$ with $($discrete$)$ asymptotic types $\P$ and $\Q$
$($associated with the weight data $(\beta,\Theta)$ and
$(-\gamma,\Theta),$ respectively$)$ if $g(y,\eta)$ has the
properties
\begin{eqnarray}\label{g}
g(y,\eta)\in
S^\mu_\clas\big(\Omega\times\R^q;\KsgX,\S^{\beta}_\P(X^\wedge)\big)
\end{eqnarray}
and
\begin{eqnarray}\label{gstar}
g^*(y,\eta)\in
S^\mu_\clas\big(\Omega\times\R^q;\K^{s,-\beta}(X^\wedge),\S^{-\gamma}_\Q(X^\wedge)\big)
\end{eqnarray}
for all $s\in\R.$ Here $g^*$ denotes the $(y,\eta)$-wise formal
adjoint with respect to the respective sesquilinear pairings
$$\K^{s,\beta}(X^\wedge)\times\K^{-s,-\beta}(X^\wedge)\rightarrow\C$$
induced by the $\K^{0,0}(X^\wedge)$ scalar product, for arbitrary
$s,\beta\in\R.$
\end{defn}
\begin{rem}
Observe that the Green operators of type $0$ in the calculus of
classical $($pseudo-differential$)$ boundary value problems are
operators with special such symbols. In this case it suffices to
replace $\KsgX$ and $\K^{s,-\beta}(X^\wedge)$ by $L^2(\R_+)$ $($with
$\R_+$ being the inner normal to the boundary in consideration$)$
and $\S^{\beta}_\P(X^\wedge)$ and $\S^{-\gamma}_\Q(X^\wedge)$ by
$\S(\Rr),$ cf. \em\cite{Schu16}.
\end{rem}
\begin{rem}
The conditions $(\ref{g})$ and $(\ref{gstar})$ are slightly stronger
than necessary. It suffices to require them for $s=0;$ however, this
is not the main point of our consideration. What we can see
immediately is that it suffices to require the conditions
$(\ref{g})$ and $(\ref{gstar})$ for all $s\in\Z$ owed by the
interpolation property of the spaces $\KsgX$ in $s.$ It follows that
the space of Green symbols of order $\mu$ and fixed $\P,\Q$ is a \Fr
space.
\end{rem}

\nt From the Green symbols which are known from the calculus of
operators on a manifold with edges we know in fact more, namely,
that the spaces $\KsgX$ and $\K^{s,-\beta}(X^\wedge)$ may even be
replaced by $\ang{r}^j\KsgX$ and $\ang{r}^j\K^{s,-\beta}(X^\wedge),$
respectively, for arbitrary $j\in\N.$ Therefore, we start with that
property. In that case it is known that the kernels of the
homogeneous components $g_{(\mu-j)}$ are $C^\infty$ functions of
$(y,\eta)\in\Omega\times(\R^q\setminus\set{0})$ with values in the
space
$$\set{\S^{\beta}_\P(X^\wedge)\TenP\S^{-\gamma}(X^\wedge)}\cap
\set{\S^\beta(X^\wedge)\TenP\S^{-\gamma}_{\cl{\Q}}(X^\wedge)}$$
where $\S^\beta(X^\wedge)=\Projlim{j\in\N}
\ang{r}^{-j}\K^{\infty,\beta}(X^\wedge)$ for any $\beta\in\R.$ Here
$\cl\Q:=\set{(\cl q_j,n_j,\cl L_j)}_j$ when
$\Q=\set{(q_j,n_j,L_j)}_j,$ and $\TenP$ denotes the (completed)
projective tensor product between the respective \Fr spaces. In this
conclusion we employ the fact that when an operator $g:H\rightarrow
F$ is continuous from a Hilbert space $H$ to a nuclear \Fr space $F$
(written as $\Projlim{j\in\N}F_j$ for Hilbert spaces $F_j$ with
nuclear embeddings $F_{j+1}\hookrightarrow F_j$ for all $j$), the
operator $g$ has a kernel in $\Projlim{j\in\N}F_j\TenH
H^*=\Projlim{j\in\N}F_j\TenP H^*=F\TenP H^*,$ cf. \cite{Jarc1}.

\section{Discrete asymptotics}

{\setlength{\baselineskip}%
{.8\baselineskip} {\Small \nt Green operators of the cone algebra
can be described by integral kernels in a specific kind of tensor
products between spaces of functions with asymptotics near $r=0$ and
of Schwartz type for $r\rightarrow\infty.$ Green symbols of the edge
calculus are symbols with values in such operators on the cone, in
this section with discrete asymptotics. We derive a characterisation
for such symbols in terms of kernels, depending on variables and
covariables on the edge, where the edge covariable appear as a
product with the axial variables $r$ and $r',$ respectively. } \par}

\subsection{Green symbols}\label{2.1}

Let $f(r,x,r',x';y,\eta)$ be a function in the space
\begin{equation}\label{X}
\set{\S^{\beta}_\P(X^\wedge)\TenP\S^{-\gamma}(X^\wedge)}\cap
\set{\S^\beta(X^\wedge)\TenP\S^{-\gamma}_{\cl{\Q}}(X^\wedge)} \TenP
S^{\mu+n+1}_\clas(\Omega\times\R^q),
\end{equation}
$\gamma,\beta,\mu\in\R,$ $n=\dim X,$ and let $\eta\mapsto[\eta]$
denote any strictly positive $C^\infty$ function in $\R^q$ such that
$[\eta]=|\eta|$ for $|\eta|>C$ for some constant $C>0.$ Form the
operator function
\begin{equation}\label{gg}
g(y,\eta)u(r,x):=\Int_X\Int_0^\infty
f(r[\eta],x,r'[\eta],x';y,\eta)u(r',x')(r')^ndr'dx'.
\end{equation}
Then we get a Green symbol of order $\mu$ in the sense of Definition
\ref{Green symbol}. For purposes below we set
\begin{equation}\label{gf}
g_f(y,\eta):=g(y,\eta).
\end{equation}
\begin{thm}\label{THM2.1.1}
Every Green symbol $g(y,\eta)$ of order $\mu$ as in
$\mathrm{Definition}$ $\ref{Green symbol}$ has a representation of
the form $(\ref{gg})$ for an element $f(r,x,r',x';y,\eta)$ in the
space $(\ref{X}).$
\end{thm}
\Pr For convenience we consider a Green symbol with constant
coefficients, i.e., $g=g(\eta)$ (the straightforward generalisation
of arguments to the $y$-dependent case will be omitted). First
observe that a simple composition of $g$ with suitable powers in $r$
and $r'$, allows us to consider the case $\beta=\gamma=0.$ Moreover,
without loss of generality we may assume $\mu=0$ (it suffices to
replace $g$ by $[\eta]^{-\mu}g$). In other words we start with $g\in
S^0_\clas\big(\R^q;\K^{s,0}(X^\wedge),\S^0_\P(X^\wedge)\big)$ with
the homogeneous components $g_{(-j)}(\eta),$ $j\in\N.$ We use the
fact that the series
\begin{equation}\label{con}
\til
g_l(\eta):=\Summa_{j=l}^\infty\chi(\frac{\eta}{c_j})g_{(-j)}(\eta)
\end{equation}
converges in
$S^{-l}_\clas\big(\R^q;\K^{s,0}(X^\wedge),\S^0_\P(X^\wedge)\big)$
for every $l\in\N.$ Here $\chi(\eta)$ is any excision function in
$\R^q,$ and $c_j$ are constants tending to $\infty$ sufficiently
fast. Then $g(\eta)-\til g(\eta),$ for $\til g(\eta):=\til
g_0(\eta),$ is of order $-\infty$ in the sense of the first part of
Definition \ref{Green symbol}. In a similar manner we can proceed
with the formal adjoint and choose, if necessary, the constants
$c_j$ once again larger, such that $g^*(\eta)-\til g^*(\eta)$ is of
order $-\infty$ in the sense of
the second part of Definition \ref{Green symbol}.\\

\nt Setting
\begin{equation}\label{S3}
\S_\P^0(X^\wedge)\TenG\S^0_{\cl{\Q}}(X^\wedge):=
\big\{\S^0_{\P}(X^\wedge)\TenP\S^0(X^\wedge)\big\}\cap
\big\{\S^0(X^\wedge)\TenP\S^0_{\cl\Q}(X^\wedge)\big\},
\end{equation}
the components $g_{(-j)}(\eta)$ can be identified with an
$\eta$-dependent kernel function of the form
$\abs{\eta}^{n+1-j}e_{(-j)}(r\abs{\eta},x,r'\abs{\eta},x';\displaystyle\frac{\eta}{\abs{\eta}}),$
$\eta\ne 0,$ for
$e_{(-j)}(r\abs{\eta},x,r'\abs{\eta},x';\displaystyle\frac{\eta}{\abs{\eta}})\in
C^\infty\big(S^{q-1}_\eta,
\S_\P^0(X^\wedge)\TenG\S^0_{\cl{\Q}}(X^\wedge)\big),$ with $S^{q-1}$
being the unit sphere in $\R^q,$ such that
\begin{equation}\label{I}
\begin{array}{rl}
  g_{(-j)}(\eta)u(r,x)=\\
  \\
  &\hspace{-2cm}\abs{\eta}^{n+1-j}\Int_X\Int_0^\infty
e_{(-j)}(r\abs{\eta},x,r'\abs{\eta},x';\frac{\eta}{\abs{\eta}})
u(r',x')(r')^ndr'dx'.
\end{array}
\end{equation}

\nt If $E$ is a \Fr space with the countable semi-norm system
$(p_k)_{k\in\N}$ we denote by $S^\mu(\R^q,E)$ the set of all $a\in
C^\infty(\R^q,E)$ such that
$$\sup_{\eta\in\R^q}{\ang\eta}^{-\mu+|\alpha|}p_k(D_\eta^\alpha a)<\infty$$
for all $\alpha\in\N^q,$ $k\in\N.$ There is then the subspace
$S_\clas^\mu(\R^q,E)$ of classical $E$-valued symbols in terms of
asymptotic expansions of elements $\chi(\eta)a_{(\mu-j)}(\eta)$ with
homogeneous components $a_{(\mu-j)}(\eta)\in
C^\infty(\R^q\setminus\{0\},E)$ of order
$\mu-j.$\\

\nt Setting
\begin{equation}\label{sj}
h_j(r,x,r',x';\eta):=\chi(\frac{\eta}{c_j})[\eta]^{n+1-j}e_{(-j)}(r,x,r',x';\frac{\eta}{\abs{\eta}})
\end{equation}
we obtain elements
\begin{equation}\label{hj}
h_j\in S_\clas^{n+1-j}\big(\R^q,
\S_\P^0(X^\wedge)\TenG\S^0_{\cl\Q}(X^\wedge)\big).
\end{equation}
Choosing the constants $c_j>0$ increasing sufficiently fast as
$j\rightarrow\infty$ we obtain convergence of
$a_l(\eta):=\Summa^\infty_{j=l}h_j(\eta)$ in
$S^{n+1-l}_\clas\big(\R^q,\S_\P^0(X^\wedge)\TenG\S^0_{\cl{\Q}}(X^\wedge)\big)$
for every $l\in\N.$ Clearly we can take the same constants as in
(\ref{con}); it suffices to take the maximums of both choices. Note
that $h_j(r,x,r',x';\eta)$ may be replaced by
$$h_j(r[\eta],x,r'[\eta],x';\eta)=\chi(\frac{\eta}{c_j})\abs\eta^{n+1-j}
e_{(-j)}(r\abs{\eta},x,r'\abs{\eta},x';\frac{\eta}{\abs{\eta}})$$
when we choose $c_0$ sufficiently large and $c_j>c_0$ for all $j\ge
1.$ According to (\ref{gf}) we obtain associated Green symbols
$g_{h_j}(\eta),$ and $\Summa^\infty_{j=l}g_{h_j}(\eta)$ converges to
$g_{a_l}(\eta)$ in the \Fr space of Green symbols of order $-l$ for
the given fixed $\P,\Q;$ this holds for every $l\in\N.$ Thus it
follows that $c(\eta):=g(\eta)-g_{a_0}(\eta)$ is a Green symbol of
order $-\infty.$ It remains to prove that there is an
$m(r,x,r',x';\eta)
\in\S\big(\R^q,\S_\P^0(X^\wedge)\TenG\S^0_{\cl{\Q}}(X^\wedge)\big)$
such that $c(\eta)=g_m(\eta).$ The Green symbol $c(\eta)$ is of
order $-\infty;$ then there is a
\begin{equation}\label{k}
k(r,x,r',x';\eta)\in\S\big(\R^q,\S_\P^0(X^\wedge)\TenG\S^0_{\cl{\Q}}(X^\wedge)\big)
\end{equation}
such that
\begin{equation*}
c(\eta)u(r,x)=\Int_X\Int_0^\infty
k(r,x,r',x';\eta)u(r',x')(r')^ndr'dx'.
\end{equation*}

\nt In Lemma \ref{k-m} below we will show that
\begin{equation}\label{m}
k\Big(\frac{r}{[\eta]},x,\frac{r'}{[\eta]},x';\eta\Big)=:m(r,x,r',x';\eta)
\in\S\big(\R^q,\S_\P^0(X^\wedge)\TenG\S^0_{\cl{\Q}}(X^\wedge)\big).
\end{equation}

\nt Then we obviously obtain $c(\eta)=g_m(\eta).$\eop
\begin{lem}\label{k-m}
We have $(\ref{k})\Rightarrow(\ref{m}).$
\end{lem}
\Pr The proof is elementary though voluminous. Therefore, we only
describe the typical steps. By virtue of (\ref{S3}) it suffices to
show that
\begin{equation}\label{K}
k\in\S\big(\R^q,\S_\P^0(X^\wedge)\TenP\S^0(X^\wedge)\big)\Rightarrow
m\in\S\big(\R^q,\S_\P^0(X^\wedge)\TenP\S^0(X^\wedge)\big)
\end{equation}
and a similar relation for Schwartz functions with values in the
second space of (\ref{S3}). Let us consider, for instance, the case
(\ref{K}). We now observe that
$$\S_\P^0(X^\wedge)\TenP\S^0(X^\wedge)=\Projlim{N\in\N}A^N\otimes_H B^N$$
for the spaces $A^N:=A^N_\P$ and $B^N,$ cf. Section \ref{1.1}, with
$\otimes_H$ being the Hilbert tensor product. Then we have
$$\S\Set{\R^q,\S_\P^0(X^\wedge)\TenP\S^0(X^\wedge)}=\Projlim{N\in\N}\S\Set{\R^q,A^N\otimes_H B^N}.$$

\nt As the semi-norm system for this space we can take
\begin{equation}\label{N}
\sup_{\eta\in\R^q}\norm{{\ang\eta}^lD^\beta_\eta
k\Big(\frac{r}{[\eta]},x,\frac{r'}{[\eta]},x';\eta\Big)}_{A^N\otimes_H
B^N}
\end{equation}
for all $l,N\in\N,$ $\beta\in\N^q.$\\

\nt It suffices to show that for every $l,\beta,N$ there are
finitely many triples $(l',\beta',N')$ such that
\begin{equation}\label{L}
\sup_{\eta\in\R^q}\norm{{\ang\eta}^{l'}D^{\beta'}_\eta
k\Big(\frac{r}{[\eta]},x,\frac{r'}{[\eta]},x';\eta\Big)}_{A^{N'}\otimes_H
B^{N'}}<\infty
\end{equation}
for all those $(l',\beta',N')$ implies that (\ref{N}) is finite.\\

\nt Let us look at the case $q=1$ and $n=\dim X=0;$ the general case
is completely analogous. For $\beta=0$ we use the fact that when
$\set\Kl_{\lambda\in\R_+}$ is a strongly continuous group of
isomorphisms on a Hilbert space $E,$ there are constants $c,M>0$
such that
\begin{equation}\label{lambda}
\norm{\Kl}_{\L(E,\til E)}\le c\Set{\max(\lambda,\lambda^{-1})}^M
\end{equation}
for all $\lambda\in\R_+.$ From (a slight modification of) Remark
\ref{kappa} we know that $u(r,x)\mapsto u(\lambda r,x),$
$\lambda\in\R_+,$ induces strongly continuous groups of isomorphisms
on the spaces $A^N$ and $B^N$ for all $N\in\N.$ Then (\ref{lambda})
yields estimates of the kind
\begin{equation}\label{est1}
\norm{k\Big(\frac{r}{[\eta]},\frac{r'}{[\eta]};\eta\Big)}_{A^N\otimes_H
B^N}\le c\ang\eta^M\norm{k(r,r';\eta)}_{A^N\otimes_H B^N}
\end{equation}
for all $\eta,$ for suitable constants $c,M>0,$ for all $k\in
A^N\otimes_H B^N.$ This gives us immediately the conclusion
$(\ref{L})\Rightarrow(\ref{N})$ with $\beta'=0$ and $l'=l+M.$\\

\nt Let us now assume $\beta=1.$ In this case we obtain
$$\frac{d}{d\eta}k\Big(\frac{r}{[\eta]},\frac{r'}{[\eta]};\eta\Big)=
\Set{\set{\varphi r\partial_r+\varphi
r'\partial_{r'}+\partial_\eta}k}\Big(\frac{r}{[\eta]},\frac{r'}{[\eta]};\eta\Big)$$
with a uniformly bounded function $\varphi(\eta)$ and
$\partial/\partial\eta$ denoting the derivative in the third
variable. Then
$$\begin{array}{l}
\displaystyle\norm{{\ang\eta}^l\Big(\frac{d}{d\eta}
k\Big)\Big(\frac{r}{[\eta]},\frac{r'}{[\eta]};\eta\Big)}_{A^N\otimes_H
B^N} \le\quad c\displaystyle\norm{{\ang\eta}^l (r\partial_{
r}k)\Big(\frac{r}{[\eta]},\frac{r'}{[\eta]};\eta\Big)}_{A^{N}\otimes_H
B^{N}}\\
\\
+\quad c\displaystyle\norm{{\ang\eta}^l(r'\partial_{
r'}k)\Big(\frac{r}{[\eta]},\frac{r'}{[\eta]};\eta\Big)}_{A^{N}\otimes_H
B^{N}}+\quad\displaystyle\norm{\ang\eta^l\Set{\partial_\eta k}
\Big(\frac{r}{[\eta]},\frac{r'}{[\eta]};\eta\Big)}_{A^{N\otimes_H
B^{N}}}
\end{array}$$
with some $c>0.$ The operator $r\partial_r$ is continuous in the
sense
\begin{equation}\label{rDr}
r\partial_r:A^N\rightarrow A^{N-1},\quad B^N\rightarrow B^{N-1}
\end{equation}
for every $N\ge1.$ In combination with the estimates (\ref{est1})
this implies

$$\begin{array}{l}
\displaystyle\norm{{\ang\eta}^l\Big(\frac{d}{d\eta}
k\Big)\Big(\frac{r}{[\eta]},\frac{r'}{[\eta]};\eta\Big)}_{A^N\otimes_H
B^N} \le\quad
c\displaystyle\norm{{\ang\eta}^{l+M}k\Set{r,r';\eta}}_{A^{N-1}\otimes_H
B^{N}}\\
\\
+\quad
c\displaystyle\norm{{\ang\eta}^{l+M}k\Set{r,r';\eta}}_{A^{N}\otimes_H
B^{N-1}}+\quad\displaystyle\norm{{\ang\eta}^{l+M}\Set{\partial_\eta
k}\Set{r,r';\eta}}_{A^{N}\otimes_H B^{N}}.
\end{array}$$

\nt The desired estimate then follows from Remark \ref{kappa} (i).
In an analogous manner we can treat the semi-norms with higher
$\eta$-derivatives.\eop\\
\begin{rem}
$\mathrm{Theorem}$ $\ref{THM2.1.1}$ remains true in analogous form
if we replace $[\eta]$ in the formula $(\ref{gg})$ by any other
strictly positive $C^\infty$ function $p(\eta)$ such that
$c[\eta]\le p(\eta)\le c' [\eta]$ for all $\eta,$ with suitable
constants $0<c<c'.$ In particular, we may take $p(\eta)=\ang\eta.$
\end{rem}

\subsection{Trace and potential symbols}\label{2.2}

The Definition \ref{Green symbol} can be generalised to $2\times2$
block matrix-valued functions $g(y,\eta)\in
S^\mu_\clas\Set{\Omega\times\R^q;\K^{0,\gamma}(X^\wedge)\oplus\C,\K^{0,\beta}(X^\wedge)\oplus\C}$
such that
$$g(y,\eta)\in
S^\mu_\clas\big(\Omega\times\R^q;\KsgX\oplus\C,\S^\beta_\P(X^\wedge)\oplus\C)$$
and
$$g^{*}(y,\eta)\in
S^\mu_\clas\big(\Omega\times\R^q;\K^{s,-\beta}(X^\wedge)\oplus\C,\S^{-\gamma}_\Q(X^\wedge)\oplus\C\big)$$
for all $s\in\R,$ with suitable $g$-dependent discrete asymptotic
types $\P,\Q$ (the pointwise adjoints refer to corresponding
sesquilinear pairings induced by the scalar product of
$\K^{0,0}(X^\wedge)\oplus\C).$ In $\C$ we always assume the trivial
group action, i.e., $\Kl$ acts as the identity for all
$\lambda\in\R_+.$ Writing
$g(y,\eta)=\Set{g_{ij}(y,\eta)}_{i,j=1,2},$ we call $g_{21}(y,\eta)$
a trace symbol and $g_{12}(y,\eta)$ a potential symbol of order
$\mu\in\R$ (of the edge calculus). Clearly $g_{22}(y,\eta)$ is
nothing other than a classical scalar symbol
of order $\mu.$\\

\nt Examples of trace and potential symbols may be obtained by
functions in
\begin{equation}\label{f21}
\S^{-\gamma}_{\cl\Q}(X^\wedge)\TenP
S^{\mu+\frac{n+1}2}_\clas(\Omega\times\R^q)\ni f_{21}(r',x';y,\eta)
\end{equation}
and
\begin{equation}\label{f12}
\S^\beta_\P(X^\wedge)\TenP
S^{\mu+\frac{n+1}2}_\clas(\Omega\times\R^q)\ni f_{12}(r,x;y,\eta),
\end{equation}
respectively; as usual, $n=\dim X.$ The symbols themselves are
obtained by integral representations of the kind
\begin{equation}\label{g21}
g_{21}(y,\eta)u=\Int_X\Int_0^\infty
f_{21}(r'[\eta],x';y,\eta)u(r',x')(r')^ndr'dx',
\end{equation}
$u(y,\eta)\in\KsgX,$ and
\begin{equation}\label{g12}
g_{12}(y,\eta)c(r,x)=cf_{12}(r[\eta],x;y,\eta),
\end{equation}
$c\in\C,$ respectively.
\begin{thm}
\begin{itemize}
\item[$\mathrm{(i)}$] Every trace symbol $g_{21}(y,\eta)$ can be
written in the form $(\ref{g21})$ for an element $(\ref{f21});$
\item[$\mathrm{(ii)}$] every potential symbol $g_{12}(y,\eta)$ can
be written in the form $(\ref{g12})$ for an element $(\ref{f12}).$
\end{itemize}
\end{thm}

\nt The proof employs analogous arguments as those for Theorem
\ref{THM2.1.1}.

\section{Continuous asymptotics}

{\setlength{\baselineskip}%
{.8\baselineskip} {\Small \nt Continuous asymptotics are motivated
by variable discrete asymptotics, where in general the asymptotic
data depend on the edge variable $y$ (i.e., exponents as well as the
occurring logarithmic powers and the coefficients may be
discontinuous in $y$). We give a definition of weighted spaces with
continuous asymptotics and write them as projective limits of
suitable Hilbert spaces. Moreover we illustrate the connection of
such variable asymptotic data with varying and branching poles of
meromorphic Melllin symbols of the calculus. Finally we give a
kernel characterisation of Green edge symbols with continuous
asymptotics similarly as in the discrete case. } \par}

\subsection{Green symbols}
In Section \ref{1.1} we have formulated spaces $\KsgXAs{\P}$ with
discrete asymptotics of type $\P$ for $r\rightarrow0,$ cf. the
formula (\ref{new}). As is known, cf. \cite{Schu2}, asymptotics of
that form can also be written as
\begin{equation}\label{u}
u(r,x)\sim\Summa_j\ang{\zeta_j,r^{-z}}
\end{equation}
where $\zeta_j$ are $C^\infty(X)$-valued analytic functionals
carried by the points $p_j\in\C$ which are of finite order (in fact,
derivatives of the Dirac distribution at $p_j$ of order
$m_j+1$ in the notation of the formula (\ref{new})).\\

\nt For an open $\U\subseteq\C$ and for a \Fr space $E,$ we denote
by $\A(\U,E)=\A(\U)\TenP E$ the space of all holomorphic $E$-valued
functions in $\U.$ Moreover, let $\A'\Set{K,C^\infty(X)}$
($=\A'(K)\TenP C^\infty(X)$) denote the space of all analytic
functionals carried by a compact set $K\subset\C.$ From generalities
on analytic functionals it follows that every
$\zeta\in\A'\Set{K,C^\infty(X)}$ can be represented in the form
\begin{equation}\label{zeta}
\zeta:h\mapsto\frac{1}{2\pi i}\Int_Cf(z)h(z)dz
\end{equation}
for some $f\in\A\Set{\C\setminus K,C^\infty(X)},$ where $C$ is a
$C^\infty$ curve counter clockwise surrounding $K$ (such that the
winding number with respect to every $z\in K$ is equal to 1). In
other words, to express (\ref{u}) it suffices to represent $\zeta_j$
by a meromorphic function with a pole at $p_j$ of order $m_j+1$ and
Laurent coefficients belonging to the space $L_j$ (cf. the
notation in Section \ref{1.1}).\\

\nt Now an element $u(r,x)\in\KsgX$ is said to have continuous
asymptotics (first in a finite weight strip $\Theta$) if there is an
element $\zeta\in\A'\Set{K,C^\infty(X)}$ for a suitable compact
$K\subset\set{z:\Re z<\frac{n+1}{2}-\gamma},$ $n=\dim X,$ such that
$$u(r,x)=\omega(r)\ang{\zeta,r^{-z}}+u_\Theta(r,x)$$
for some $u_\Theta\in\KsgXAs\Theta;$ here, as usual, $\omega$ is a
cut-off function.\\

\nt In order to unify notation in connection with discrete or
continuous asymptotics we consider the space
\begin{equation}\label{EP}
\set{\omega(r)\ang{\zeta,r^{-z}}:\zeta\in\A'(K,C^\infty(X))}.
\end{equation}

\nt The quotient space of (\ref{EP}) with respect to the equivalence
relation $u\sim v\Leftrightarrow u-v\in\KsgXAs\Theta$ is called a
continuous asymptotic type $\P,$ associated with weight data
$(\gamma,\Theta).$ The cut-off function $\omega$ is fixed, but the
quotient space is independent of $\omega.$ Denoting the space
(\ref{EP}) by $\E_\P(X^\wedge)$ we then define
\begin{equation}\label{nds}
\KsgXAs\P:=\KsgXAs\Theta+\E_\P(X^\wedge)
\end{equation}
in the \Fr topology of the non-direct sum. To recall the
terminology, the non-direct sum of two \Fr spaces $E$ and $F$
(embedded in a Hausdorff topological vector space) is defined as
$E+F:=\set{e+f:e\in E,\;f\in F}$ endowed with the \Fr topology from
$E+F\simeq E\oplus F/\Delta,$ where $\Delta:=\set{(e,-e):e\in E\cap
F}.$ Note that $\P$ only depends on the set $K\cap\set{\Re
z>\frac{n+1}{2}-\gamma+\vartheta}.$ For a generalisation to infinite
weight intervals we define the system $\VV$ of closed subsets
$V\subset\C$ such that $V\cap\set{c\le\Re z\le c'}$ is compact for
every $c\le c'.$ Then for every $V\in\VV$ contained in $\set{\Re
z<\frac{n+1}{2}-\gamma}$ we can consider $V_\vartheta:=V\cap\set{\Re
z\ge\frac{n+1}{2}-\gamma+\vartheta-1}$ and the associated continuous
asymptotic type $\P_\vartheta.$ We then have continuous embedding
$$\KsgXAs{\P_\vartheta}\hookrightarrow\KsgXAs{\P_{\vartheta'}}$$
for every $\vartheta<\vartheta',$ and we then set
$\KsgXAs\P:=\displaystyle\lim_{\substack\longleftarrow}\KsgXAs{\P_\vartheta}$
in the \Fr topo-logy of the projective limit taken over any
monotonically decreasing sequence of negative reals $\vartheta.$ The
subscript $\P$ incorporates a continuous asymptotic type associated
with $(\gamma,\Theta)$ for $\Theta=(-\infty,0],$ and stands for the
equivalence class represented by such a sequence
$\set{\P_\vartheta}$ with $\vartheta$ running over any monotone
sequence tending to $-\infty.$ The equivalence relation just means
the equality of the respective projective limits. The set $V$ is
called a carrier set of the asymptotic type $\P$ (when
$\vartheta=-\infty,$ otherwise $V\cap\set{\Re
z\ge\frac{n+1}2-\gamma+\vartheta}$ is
called the carrier of the corresponding $\P).$\\

\nt We do not need the sets $V$ in full generality. Let us content
ourselves with those $V$ that are convex in imaginary direction,
i.e., $z_0,z_1\in V$ and $\Re z_0=\Re z_1$ imply $\lambda
z_0+(1-\lambda)z_1\in V$ for all $0\le\lambda\le1.$ There is then an
obvious one-to-one correspondence between such $V$ contained in
$\set{\Re z<\frac{n+1}{2}-\gamma}$ and associated continuous
asymptotic types by the above construction.\\

\nt If $\P$ is a continuous asymptotic type, we set
$$\S^\gamma_\P(X^\wedge):=\Projlim{N\in\N}\ang{r}^{-N}\K_\P^{\infty,\gamma}(X^\wedge)$$
which is a nuclear \Fr space in the topology of the projective limit.\\

\nt Let us make some remarks about the motivation of continuous
asymptotics. As noted in the introduction the elliptic regularity of
solutions to elliptic equations $Au=f$ on a wedge
$X^\wedge\times\Omega,$ $\Omega\subseteq\R^q$ open, and $A$
edge-degenerate of the form (\ref{neu}), contains a statement on
asymptotics of $u(r,x,y)$ for $r\rightarrow0,$ even if we are
considering $C^\infty$ functions on $X^\wedge\times\Omega.$
Similarly as (\ref{new}) the asymptotics have the form
\begin{equation}\label{uu}
u(r,x,y)\sim\Summa_j\Summa^{m_j(y)}_{k=0}c_{jk}(x,y)r^{-p_j(y)}\log^kr
\qquad\mathrm{for}\quad r\rightarrow0,
\end{equation}
where the exponents $-p_j(y)$ and the numbers $m_j(y)$ are
determined by those points $z\in\C$ where the operators
(\ref{sigmac}) are not bijective, cf. \cite{Schu28}. These points
(as well as the $m_j(y)$) may depend on $y$ in a very irregular way.
This may happen even for $n=\dim X=0.$ The inverse of (\ref{sigmac})
is then a family of meromorphic functions, and the main ingredients
of the parametrices $P$ of $A$ are Mellin operators with such
symbols. Applying $P$ to functions (say, with compact support with
respect to $r\in\R_+$) gives us functions $u(r,x,y)$ of a behaviour
like (\ref{uu}). If we consider the Mellin transform $\Set{M(\omega
u)}(z,x,y)$ (for any cut-off function $\omega(r)$ on the half-axis)
we obtain a family of meromorphic functions in the complex plane the
poles and multiplicities of which inherit the corresponding
behaviour of the Mellin symbols. Now the parameter-dependent
asymptotics of $u$ can be interpreted in terms of functions
$\zeta(y)\in C^\infty\Set{\Omega,\A'(K,C^\infty(X))}$ for suitable
compact $K,$ such that $\zeta(y)$ is pointwise discrete and of
finite order but of the above mentioned irregular behaviour. Here
\virg{pointwise discrete} means that $\ang{\zeta(y),r^{-z}}$ has the
form (\ref{uu}) for certain $p_j\in K,$ $m_j\in\N$ for every
$y\in\Omega.$
\begin{prop}\label{3.1.1}
For every continuous asymptotic type $\P$ associated with weight
data $(\gamma,\Theta)$ there is a scale of Hilbert spaces $A^N_\P,$
$N\in\N,$ with nuclear embeddings $A^N_\P\hookrightarrow A^{N-1}_\P$
for every $N\ge1$ such that
$$\S^\gamma_\P(X^\wedge)=\Projlim{N\in\N}A^N_\P.$$
The spaces $A^N_\P$ can be chosen as continuously embedded subspaces
of $\K^{0,\gamma}(X^\wedge)$ such that $(\ref{kl})$ induces a
strongly continuous group of isomorphisms
$$\Kl:A^N_\P\rightarrow A^N_\P$$
for every $N\in\N.$
\end{prop}

\Pr We first assume the weight interval $\Theta$ to be finite.
Similarly as (\ref{nds}) we write $\S_\P^\gamma(X^\wedge)$ as a
(non-direct) sum of \Fr spaces, namely,
$$\S_\P^\gamma(X^\wedge)=\Projlim{N\in\N}\ang{r}^{-N}\K_\Theta^{\infty,\gamma}(X^\wedge)+\E_\P(X^\wedge).$$
We then consider the spaces
\begin{equation}\label{BNcont}
\til
A_\P^N:=\ang{r}^{-N}\K^{N,\gamma-\vartheta-\frac1{N+1}}(X^\wedge)+\E_\P^N(X^\wedge).
\end{equation}
The meaning of the first summand is clear, cf. also the formula
(\ref{KTheta}); so it remains to define $\E_\P^N(X^\wedge).$ Recall
that the space (\ref{EP}) may be described in terms of a compact set
$K:=V\cap\set{\Re z\ge\frac{n+1}2-\gamma+\vartheta-1}$ for a set
$V\subset\set{\Re z<\frac{n+1}2-\gamma}$ of the above mentioned
kind. Choose any $C^\infty$ curve $C_N\subset\C\setminus K$ counter
clockwise surrounding $K,$ such that the winding number with respect
to any $z\in K$ is equal to $1$ and $\mathrm{dist}(z,K)\le\frac1N$
(it is well known that such curves always exist). Then the weighted
Mellin transforms $\Set{M_{\gamma-\frac n2}u}|_{C_N}$ belong to
$C^\infty(C_N)$ (concerning notation about the weighted Mellin
transform, cf. Section \ref{1.01}). By $\E_\P^N(X^\wedge)$ we then
denote the completion of $\set{\Set{M_{\gamma-\frac
n2}u}|_{C_N}:u\in\E_\P(X^\wedge)}$ in the norm of $H^N(C_N),$ the
Sobolev space of smoothness $N$ on the curve $C_N.$ This is a
Hilbert space, and we have
$\E_\P(X^\wedge)=\Projlim{N\in\N}\E_\P^N(X^\wedge).$\\

\nt It is now clear that the space (\ref{BNcont}) is nuclearly
embedded into a corresponding space of analogous structure belonging
to $N'<N$ when $N-N'$ is sufficiently large. This allows us to find
a sequence of $N_j,$ $j\in\N,$ with $N_{j+1}>N_j$ such that, if we
set $A_\P^j:=\til A_\P^{N_j},$ we have nuclearity of
$A_\P^{j+1}\hookrightarrow A_\P^j$ for all $j.$ We then have
$$\S_\P^\gamma(X^\wedge)=\Projlim{j\in\N}A_\P^j.$$
It remains to note that the group action
$\Kl:u(r,.)\mapsto\lambda^\frac{n+1}2u(\lambda r,.)$ restricts to a
group action on the space $A_\P^j$ for every $j.$ It suffices to
check that for the spaces (\ref{BNcont}). The factor
$\lambda^\frac{n+1}2$ is not essential, so we have to look at the
influence of rescaling to the space $\E_\P^N(X^\wedge).$ By
definition we restrict the Mellin transform to the curve $C_N$ and
measure the result in $H^N(C_N).$ The Mellin transform of the
rescaled function is obtained by multiplying the original one by
$\lambda^{-z}.$ The continuous dependence of the $H^N(C_N)$-norm on
$\lambda\in\R_+$ is then obvious. For the infinite weight interval
$\Theta=(-\infty,0]$ we first write
$\S^\gamma_\P(X^\wedge)=\Projlim{m\in\N}\S^\gamma_{\P_{\vartheta_m}}(X^\wedge)$
for a sequence of finite $\vartheta_m<0$ tending to $-\infty$ and
form the spaces $A^N_{\P_{\vartheta_m}}$ for every $\vartheta$ of
this sequence, such that
$\S^\gamma_{\P_{\vartheta_M}}=\Projlim{N\in\N}A^N_{\P_{\vartheta_M}}.$
Then we can set $A^N_\P:=A^N_{\P_{\vartheta_N}}.$ \eop
\begin{rem}\label{3.1.2}
$\mathrm{Definition}$ $\ref{Green symbol}$ has an immediate
generalisation to Green symbols with continuous asymptotic types
$\P$ and $\Q,$ associated with weight data $(\beta,\Theta)$ and
$(-\gamma,\Theta),$ respectively.
\end{rem}

\subsection{Mellin symbols}\label{3.2}
As noted in the beginning Green operators on a manifold with conical
singularities belong to the algebra of cone pseudo-differential
operators. Technically they appear as remainders in some typical
operations with so called smoothing Mellin operators, cf.
(\ref{MellinOp}), also defined in terms of asymptotic data. Recall
that $M$ denotes the Mellin transform on $\R_+,$ i.e.,
$Mu(z):=\Int_0^\infty r^{z-1}u(r)dr$ and Mellin operators on an
infinite stretched cone $X^\wedge$ occur with operator-valued
symbols taking values in the calculus of operators on the cone. In
this connection the Mellin symbols depend on edge variables and
covariables, and the mapping properties refer to asymptotic data for
$r\rightarrow0.$ It is typical that the Mellin amplitude functions
are not only defined on $\Gamma_{\frac12-\gamma}$ but in the complex
$z$-plane, up to a subset $V$ which encodes asymptotic properties,
similarly as in the context of functions with (discrete or
continuous) asymptotics. We want to give a definition and then
observe the way how Green operators are induced by Mellin
operators with asymptotics.\\

\nt From now on, we assume the sets $V\in\VV$ to be convex in
imaginary direction. A $V$-excision function is any $\chi\in
C^\infty(\C)$ such that $\chi(z)=0$ when
$\mathrm{dist}(z,V)<\eps_0,$ $\chi(z)=1$ for
$\mathrm{dist}(z,V)>\eps_1$ for certain
$0<\eps_0<\eps_1.$\\

\nt By $M_V^{-\infty}(X)$ we denote the space of all
$f(z)\in\A\Set{\C\setminus V,L^{-\infty}(X)}$ such that
$$\chi(z)f(z)|_{\Gamma_\beta}\in\S\Set{\Gamma_\beta,L^{-\infty}(X)}$$
for every $V$-excision function $\chi$ and every real $\beta,$
uniformly in compact $\beta$-intervals. Moreover, let $M_\O^\mu(X),$
$\mu\in\R,$ denote the space of all
$h(z)\in\A\Set{\C,L^\mu_\clas(X)}$ such that
$$h(z)|_{\Gamma_\beta}\in L^\mu_\clas(X;\Gamma_\beta)$$
for every real $\beta,$ uniformly in compact $\beta$-intervals.\\

\nt The spaces $M_V^{-\infty}(X)$ and $M_\O^\mu(X)$ are nuclear
\Fr spaces in a natural way.\\

\nt Let us set
$$M_V^\mu(X):=M_\O^\mu(X)+M_V^{-\infty}(X)$$

\nt in the \Fr topology of the non-direct sum. Then, for every
$f(r,r',z)$ belonging to the space
$C^\infty\Set{\Rr\times\Rr,M^\mu_V(X)},$ we can form associated
weighted Mellin operators $\op_M^\beta(f),$ for every weight
$\beta\in\R$ such that $V\cap\Gamma_{\frac12-\beta}=\emptyset.$
\begin{thm}\label{Mellin}
For every $f(r,r',z)\in C^\infty\Set{\Rr\times\Rr,M^\mu_V(X)},$
$V\cap\Gamma_{\frac{n+1}2-\gamma}=\emptyset,$ $n=\dim X,$ the
operator $\omega\op_M^{\gamma-\frac n2}(f)\til\omega$ $($with
cut-off functions $\omega,\til\omega)$ induces continuous operators
$$\omega\op_M^{\gamma-\frac{n}2}(f)\til\omega:\KsgX\rightarrow\K^{s-\mu,\gamma}(X^\wedge)$$
and
$$\omega\op_M^{\gamma-\frac{n}2}(f)\til\omega:\KsgXAs\P\rightarrow\K^{s-\mu,\gamma}_\Q(X^\wedge)$$
for every $s\in\R$ and every continuous asymptotic type $\P$ with
some resulting continuous asymptotic type $\Q,$ associated with the
weight data $(\gamma,\Theta)$ for every $\Theta=(\vartheta,0],$
$-\infty\le\vartheta<0.$
\end{thm}

\nt This result is known, cf. \cite{Schu2}. Recall that the main
idea of the continuity in spaces with continuous asymptotics is to
characterise the Mellin transforms of $\omega u$ as holomorphic
functions outside the union of $V$ and the carrier  set of the
asymptotic type $\P;$ then we obtain another carrier set which
just determines the asymptotic type $\Q.$\\

\nt Mellin operators as in Theorem \ref{Mellin} belong to the
ingredients of parametrices of elliptic (pseudo-differential)
operators on manifolds with conical singularities (modelled on
$X^\wedge$) or edges (modelled on
$X^\wedge\hspace{-.1cm}\times\Omega$ for some open set
$\Omega\subseteq\R^q$). These operators are combined with
other operators of the calculus.\\

\nt To see how Green operators appear in the set-up of conical
singularities, consider the following (known) result; for
completeness we give a proof here.
\begin{prop}
Consider an element $f\in M^{-\infty}_V(X);$ let $j>0,$ $\mu\in\R,$
and let $\gamma-j\le\beta,\delta\le\gamma$ for some reals
$\beta,\delta,$ such that
$$V\cap\Gamma_{\frac{n+1}2-\beta}=V\cap\Gamma_{\frac{n+1}2-\delta}=
\emptyset.$$
For the operator
$$g:=\omega r^{-\mu+j}\op_M^{\beta-\frac{n}2}(f)\til\omega-\omega
r^{-\mu+j}\op_M^{\delta-\frac{n}2}(f)\til\omega,$$ $n=\dim X,$ there
are continuous asymptotic types $\P$ and $\Q$ associated with the
weight data $(\gamma-\mu,\Theta)$ and $(-\gamma,\Theta),$
respectively, such that
\begin{equation}\label{GG}
g:\ang{r}^l\KsgX\rightarrow\S_\P^{\gamma-\mu}(X^\wedge),
\end{equation}
and
\begin{equation}\label{GG*}
g^*:\ang{r}^l\K^{s,-\gamma+\mu}(X^\wedge)\rightarrow\S_\Q^{-\gamma}(X^\wedge)
\end{equation}
are continuous operators for all $s\in\R,$ $l\in\N.$
\end{prop}

\Pr Let us check the mapping property (\ref{GG}); the property
(\ref{GG*}) can be verified in an analogous manner by passing to
formal adjoints of the involved Mellin operators, using the fact
that they are of analogous type with resulting `adjoint' Mellin
symbols, etc., cf. \cite{Schu2}. By virtue of the fact that the
operators contain cut-off functions we immediately see that the
factors $\ang{r}^l$ are harmless; so we may look at the case $l=0.$
In addition it suffices to assume $\mu=0.$ The operator $g$ is then
continuous as a map $\KsgX\rightarrow\KsgX$ because of the assumed
weight conditions. We have
\begin{equation}\label{dif}
gu(r)=\frac1{2\pi
i}r^j\Int_{\Gamma_{\frac{n+1}2-\beta}}\hspace{-.8cm}r^{-z}f(z)M\til\omega
u(z)dz-\frac1{2\pi
i}r^j\Int_{\Gamma_{\frac{n+1}2-\delta}}\hspace{-.8cm}r^{-z}f(z)M\til\omega
u(z)dz.
\end{equation}
Let, for instance, $\beta\le\delta.$ Observe that $M\til\omega u(z)$
for $u\in\KsgX$ is holomorphic in $\Re z\ge\frac{n+1}2-\gamma.$
Thus, because of the position of ${\Gamma_{\frac{n+1}2-\beta}}$ and
${\Gamma_{\frac{n+1}2-\delta}}$ on the right of
${\Gamma_{\frac{n+1}2-\gamma}}$ we can replace the difference of
integrals (\ref{dif}) as an integration over a closed curve $C$
counter clockwise surrounding the compact set
$K:=V\cap\set{\frac{n+1}2-\beta<\Re z<\frac{n+1}2-\delta}.$ The
function $\til f(z):=f(z)M\til\omega u(z)$ is holomorphic in the
strip outside $K.$ Hence (\ref{dif}) takes the form (\ref{zeta}) for
$h(z)=r^{-z},$ up to the factor $\omega(r)r^j.$ We thus obtain
altogether $gu(r)=\omega(r)r^j\ang{\zeta,r^{-z}}$ for
$\zeta\in\A'(K,C^\infty(X))$ which gives us the mapping property
(\ref{GG}), where the asymptotic type $\P$ is represented by the
compact set $K,$ cf. the notation in connection with (\ref{EP}).
\eop\\

\nt Let us consider what are called Mellin edge symbols. Such
symbols are finite linear combinations of operator families of the
form
\begin{equation}\label{Me}
m(y,\eta):=\omega(r[\eta])r^{-\mu+j}\op_M^{\gamma_j-\frac
n2}(f_{j\alpha})(y)\eta^\alpha\til\omega(r[\eta])
\end{equation}
for cut-off functions $\omega,\til\omega,$ and $f_{j\alpha}(y)\in
C^\infty\Set{\Omega,M_V^{-\infty}(X)}$ for a set $V\in\VV$ such that
$V\cap\Gamma_{\frac{n+1}2-\gamma_j}=\emptyset,$
$\Omega\subseteq\R^q$ open. In such expressions we have $j\in\N,$
$\alpha\in\N^q,$ $|\alpha|\le j,$ and the weights $\gamma_j\in\R$
are assumed to satisfy the condition
\begin{equation}\label{gamma}
\gamma-j\le\gamma_j\le\gamma
\end{equation}
for every $j\in\N.$ Then (\ref{Me}) is a $C^\infty$ family of
continuous operators
$$m(y,\eta):\KsgX\rightarrow\S^{\gamma-\mu}(X^\wedge),$$
cf. Section \ref{1.1}. We have, in fact, more, namely
$$m(y,\eta)\in S_\clas^\mu\Set{\Omega\times\R^q;\KsgX,\S^{\gamma-\mu}(X^\wedge)}$$
for every $s\in\R,$ cf. notation (\ref{cl}), and
$$m(y,\eta)\in S_\clas^\mu\Set{\Omega\times\R^q;\KsgXAs\P,\S^{\gamma-\mu}_\Q(X^\wedge)}$$
for every continuous asymptotic type $\P$ with some resulting
continuous asymptotic type $\Q$ (associated with the weight data
$(\gamma,\Theta)$ and $(\gamma-\mu,\Theta),$ respectively).
Moreover, the pointwise formal adjoint $m^*(y,\eta)$ (cf. also
Definition \ref{Green symbol}), has a similar structure as
(\ref{Me}), i.e., we have
$$m^*(y,\eta)\in S_\clas^\mu\big(\Omega\times\R^q;\K_{(\RR)}^{s,-\gamma+\mu},\S^{-\gamma}_{(\S)}(X^\wedge)\big)$$
for all $s\in\R,$ where the subscripts mean without or with the
corresponding continuous asymptotic types.\\

\nt There are now several essential operations in the edge symbolic
calculus which produce Green symbols in the sense of Remark
\ref{3.1.2}. More precisely, we obtain Green symbols $g(y,\eta)$ of
the kind
\begin{equation}\label{gcont}
g(y,\eta)\in
S^\mu_\clas\Set{\Omega\times\R^q;\ang{r}^l\KsgX,\S_\P^{\gamma-\mu}(X^\wedge)},
\end{equation}
such that
\begin{equation}\label{gstarcont}
g^*(y,\eta)\in
S^\mu_\clas\Set{\Omega\times\R^q;\ang{r}^l\K^{s,-\gamma+\mu}(X^\wedge),\S_\Q^{-\gamma}(X^\wedge)}
\end{equation}
for suitable continuous asymptotic types $\P$ and $\Q$ associated
with the weight data $(\gamma-\mu,\Theta)$ and $(-\gamma,\Theta),$
respectively, for all $l\in\N,$ $s\in\R.$

\begin{rem} An element $m(y,\eta)$ of the form $(\ref{Me})$ is a
Green symbol for every $j>-\vartheta,$ where $\Theta=(\vartheta,0]$
is the finite weight strip which plays the role in the continuous
analogue of $\mathrm{Definition}$ $\ref{Green symbol}$.
\end{rem}

\nt Another point concerns the fact that there may be different
choices of $\gamma_j$ (when $j>0,$ otherwise for $j=0$ we have
$\gamma_0=\gamma$) such that (\ref{gamma}) holds. Let $\til\gamma_j$
denote any other choice. Then we have the following result:
\begin{rem} Let $j>0,$ and let $\gamma_j$ and $\til\gamma_j$
denote different weights satisfying $(\ref{gamma}).$ Consider the
operators $(\ref{Me})$ for both weights, e.g.,
\begin{equation*}
\til m(y,\eta):=\omega(r[\eta])r^{-\mu+j}\op_M^{\til\gamma_j-\frac
n2}(f_{j\alpha})(y)\eta^\alpha\til\omega(r[\eta]).
\end{equation*}

\nt Then we have $m(y,\eta)=\til m(y,\eta)$ modulo a Green symbol
with continuous asymptotics, cf. $\mathrm{Remark}$ $\ref{3.1.2}.$
\end{rem}

\subsection{Integral representations}

Let $f(r,x,r',x';y,\eta)$ be a function in the space
\begin{equation}\label{Xcont}
\set{\S^{\gamma-\mu}_\P(X^\wedge)\TenP\S^{-\gamma}(X^\wedge)}\cap
\set{\S^{\gamma-\mu}(X^\wedge)\TenP\S^{-\gamma}_{\cl{\Q}}(X^\wedge)}
\TenP S^{\mu+n+1}_\clas(\Omega\times\R^q),
\end{equation}
now for continuous asymptotic types $\P$ and $\Q.$ Then the integral
representation
\begin{equation}\label{ggcont}
g(y,\eta)u(r,x):=\Int_X\Int_0^\infty
f(r[\eta],x,r'[\eta],x';y,\eta)u(r',x')(r')^ndr'dx',
\end{equation}
$n=\dim X,$ gives us special Green symbols with the properties
(\ref{gcont}) and (\ref{gstarcont}) for all $l\in\N,$ $s\in\R.$
\begin{thm}\label{3.3.1}
Let $g(y,\eta)$ satisfy the conditions $(\ref{gcont})$ and
$(\ref{gstarcont})$ for all $l\in\N,$ $s\in\R.$ Then there is an
$f(r,x,r',x';y,\eta)$ in the space $(\ref{Xcont})$ such that the
integral representation $(\ref{ggcont})$ holds.
\end{thm}
\Pr The proof employs analogous steps as that of Theorem
\ref{THM2.1.1}; so we only discuss the main ideas. For simplicity we
omit again the $y$-variable and write the Green symbol $g(\eta)$ as
an asymptotic sum of the kind (\ref{con}), modulo a Green symbol of
order $-\infty.$ For the homogeneous components we take the integral
representation (\ref{I}) for all $j,$ then form the functions
(\ref{sj}) and obtain the symbols (\ref{hj}). This yields the
corresponding analogue of $g_{a_0}(\eta)$ which is of the desired
integral form, modulo a Green symbol of order $-\infty,$ given by a
kernel like (\ref{k}). It then remains to show the analogue of Lemma
\ref{k-m} for the case with continuous asymptotics. The proof of
that is a purely technical (but elementary) construction in terms of
the scales of spaces $A^N$ and $B^N.$ The spaces $B^N$ are the same
as before, while the $A^N$ are constructed in Proposition
\ref{3.1.1}. The main new aspect to be employed in the proof is the
first of the relations (\ref{rDr}). In the present case we have to
look at (\ref{BNcont}). The first summand is as in (\ref{APN}), and
it remains to observe that $-r\partial_ r$ transforms the space
$\E_\P^N(X^\wedge)$ to $\E_\P^{N-1}(X^\wedge),$ modulo a flat
contribution which is absorbed by the first summand in
(\ref{BNcont}). Applying $-r\partial_ r$ to the second factor of
$\omega(r)\ang{\zeta,r^{-z}}\in\E^N_\P(X^\wedge),$ cf. (\ref{EP}),
we obtain $\omega(r)\ang{\zeta,zr^{-z}};$ thus we remain in the
space $\E_\P(X^\wedge)$ and hence, from the continuity of
$\E_\P(X^\wedge)\rightarrow\E_\P(X^\wedge),$
$\omega(r)\ang{\zeta,r^{-z}}\mapsto\omega(r)\ang{\zeta,zr^{-z}}$ and
the definition of $\E^N_\P(X^\wedge),$ we immediately obtain the
desired relation, i.e., $-r\partial_
r:\E^N_\P(X^\wedge)\rightarrow\E^{N-1}_\P(X^\wedge).$ The other
element of the proof are very close to the ones of Lemma \ref{k-m}
and will be omitted.\eop\\

\nt Analogously as the discrete case, cf. Section \ref{2.2}, we can
consider $2\times2$ block matrix-valued functions $g(y,\eta)\in
S^\mu_\clas\Set{\Omega\times\R^q;\K^{0,\gamma}(X^\wedge)\oplus\C,\K^{0,\beta}(X^\wedge)\oplus\C}$
such that
\begin{equation}\label{MP1}
g(y,\eta)\in
S^\mu_\clas\big(\Omega\times\R^q;\K^{s,\gamma}(X^\wedge)\oplus\C,\S^\beta_\P(X^\wedge)\oplus\C\big)
\end{equation}
and
\begin{equation}\label{MP2}
g^{*}(y,\eta)\in
S^\mu_\clas\big(\Omega\times\R^q;\K^{s,-\beta}(X^\wedge)\oplus\C,\S^{-\gamma}_\Q(X^\wedge)\oplus\C\big)
\end{equation}
for all $s\in\R,$ with suitable $g$-dependent continuous asymptotic
types $\P,\Q.$ Let $g(y,\eta)=(g_{ij}(y,\eta))_{i,j=1,2};$ then we
call $g_{21}(y,\eta)$ a trace symbol and $g_{12}(y,\eta)$ a
potential symbol of order $\mu\in\R,$ while $g_{22}(y,\eta)$ is
nothing
other than a classical scalar symbol (of order $\mu$).\\

\nt Let $f_{21}\in\S^{-\gamma}_{\cl\Q}(X^\wedge)\TenP
S^{\mu+\frac{n+1}2}_\clas(\Omega\times\R^q)$ and
$f_{12}\in\S^\beta_\P(X^\wedge)\TenP
S^{\mu+\frac{n+1}2}_\clas(\Omega\times\R^q),$ and consider the
integral representations
\begin{equation}\label{g21cont}
g_{21}(y,\eta)u=\Int_X\Int_0^\infty
f_{21}(r'[\eta],x';y,\eta)u(r',x')(r')^ndr'dx',
\end{equation}
$u(r,x)\in\KsgX,$ and
\begin{equation}\label{g12cont}
g_{12}(y,\eta)c(r,x)=cf_{12}(r[\eta],x;y,\eta),
\end{equation}
$c\in\C.$ Then we have (trace and potential) symbols satisfying the
mapping properties (\ref{MP1}) and (\ref{MP2}).
\begin{thm}\label{3.3.2}
\begin{itemize}
\item[$\mathrm{(i)}$] Every trace symbol $g_{21}(y,\eta)$ can be
written in the form $(\ref{g21cont})$ for an element
$f_{21}(r',x';y,\eta)\in\S^{-\gamma}_{\cl\Q}(X^\wedge)\TenP
S^{\mu+\frac{n+1}2}_\clas(\Omega\times\R^q);$
\item[$\mathrm{(ii)}$] every potential symbol $g_{12}(y,\eta)$ can
be written in the form $(\ref{g12cont})$ for an element
$f_{12}(r,x;y,\eta)\in\S^\beta_\P(X^\wedge)\TenP
S^{\mu+\frac{n+1}2}_\clas(\Omega\times\R^q).$
\end{itemize}
\end{thm}

\nt This can be proved by analogous arguments as for Theorem
\ref{3.3.1}. Here, as usual, $n=\dim X.$

\begin{rem}
\begin{itemize}
\item[$(i)$] Green symbols in the sense of block matrices
$(\ref{MP1})$ can be composed within the respective spaces of Green
symbols $($with discrete or continuous asymptotics$)$ and the
homogeneous principal components behave multiplicatively.
\item[$(ii)$] The classes of Green symbols $($with discrete or
continuous asymptotics$)$ are closed under asymptotic summation when
the involved asymptotic types are the same for all the summands.
\end{itemize}
\end{rem}

\nt Let us conclude this section with a few intuitive remarks on the
nature of continuous asymptotics which give rise to some
`unexpected' examples of Green, trace, or potential operators in
that context. First note that when $\omega(r)$ is any fixed cut-off
function and $z\in\C,$ $\Re z<\frac{n+1}2-\gamma,$ then we have
$\omega(r)r^{-z}c(x)\in\K^{\infty,\gamma}(X^\wedge)$ for any $c\in
C^\infty(X).$ Recall that such functions may be interpreted as
singular functions of the discrete cone asymptotics for
$r\rightarrow0$ with exponent $-z.$ Now if $\zeta(y)\in
C^\infty\big(\Omega,\A'(K,C^\infty(X))\big)$ is any family of
analytic functionals carried by a compact set $K\subset\set{\Re
z<\frac{n+1}2-\gamma},$ the function
\begin{equation}\label{as}
\omega(r)\ang{\zeta(y),r^{-z}}
\end{equation}
may be interpreted as the linear superposition of singular functions
for the discrete asymptotics with the ($y$-dependent) density
$\zeta(y).$ Such `densities' may be organised as follows. Choose an
arbitrary function $f(y,z)\in C^\infty(\Omega\times(\C\setminus
K),C^\infty(X)),$ holomorphic in $z\in\C\setminus K$ that extends
for every $y\in\Omega$ to a certain ($C^\infty(X)$-valued)
meromorphic function in $z\in\C$ with poles $p_j(y)\in K,$
$j=1,\dots,N(y),$ of multiplicities $m_j(y)+1.$ These poles
including multiplicities may be not constant in $y.$ Then, setting
$$\zeta(y):h(z)\mapsto\frac1{2\pi i}\Int_Ch(z)f(y,z)dz$$
where $C\subset\set{\Re z<\frac{n+1}2-\gamma}$ is a compact (say,
$C^\infty$) curve counter clockwise surrounding the set $K$ such
that $C$ has the winding number $1$ with respect to every point of
$K,$ we obtain an element $\zeta(y)\in
C^\infty\big(\Omega,\A'(K,C^\infty(X))\big)$ such that (\ref{as})
has discrete asymptotics of the kind (\ref{uu}) for every fixed $y.$
Now
$$\C\ni c\mapsto c\cdot\omega(r[\eta])[\eta]^{\mu+\frac{n+1}2}\ang{\zeta(y),(r[\eta])^{-z}}$$
defines a potential symbol of order $\mu$ in the frame of continuous
asymptotics which just produces functions with pointwise (in $y$)
discrete but branching asymptotics. In a similar manner we can
organise trace symbols which reflect such asymptotics as well as
more general Green symbols of that kind. Constructions of that kind
may also be found in \cite{Schu2}.

\section{Green operators}

{\setlength{\baselineskip}%
{.8\baselineskip} {\Small \nt We give an idea on how the Green
operators of the edge calculus are organised, and we then have a
look at a parameter-dependent variant and formulate a result on
kernel cut-off which yields holomorphic dependence on parameters,
here for the case of constant discrete or continuous asymptotics.
This is done on the level of kernel functions constructed in the
preceding sections.} \par}

\subsection{Green operators on a manifold with edges}

\nt Let $M$ be a compact manifold with edge $Y,$ locally near any
$y\in Y$ modelled on $X^\triangle\times\R^q,$ where $X$ is a closed
compact $C^\infty$ manifold. Recall that transition functions
between (open) stretched wedges $\R_+\times X\times\R^q\ni(r,x,y)$
are assumed to be $C^\infty$ up to $r=0.$ In addition we choose the
global atlas by such singular charts near $Y$ in such a way that the
transition functions are constant with respect to $r$ for $0<r<\eps$
for some $\eps>0.$ By $\M$ we denote the stretched manifold
associated with $M,$ see the introduction.\\

\nt We consider the weighted edge Sobolev space $\WW^{s,\gamma}(\M)$
that is defined as the subspace of all $u\in
H^s_\loc(\mathrm{int}\M)$ which locally near $Y$ in the coordinates
$(r,x,y)$ belong to $\WW^s(\R^q,\KsgX).$ Here $\WW^s(\R^q,E)$ for a
Hilbert space $E$ with group action $\Kl$ is the completion of
$\S(\R^q,E)$ with respect to the norm
$$\set{\Int\ang\eta^{2s}\norm{\kappa_{\ang\eta}^{-1}\hat u(\eta)}_E^2d\eta}^\frac12,$$
with $\hat u(\eta)$ being the Fourier transform of $u$ in $\R^q.$ In
a similar manner we define $\WW^s(\R^q,E)$ for a \Fr space $E$ which
is the projective limit of Hilbert spaces $E^j$ with group actions,
with continuous embeddings $\dots\hookrightarrow
E^{j+1}\hookrightarrow E^j\hookrightarrow\cdots\hookrightarrow E^0$
for all $j\in\N,$ such that the group action on $E^j$ is the
restriction of the one on $E^0$ for every $j.$ This allows us to
define subspaces
\begin{equation}\label{WP}
\WW^s(\R^q,\S^\gamma_\P(X^\wedge))
\end{equation}
of $\WW^s(\R^q,\KsgX)$ for any (discrete or continuous) asymptotic
type $\P,$ using the fact that $\S^\gamma_\P(X^\wedge)$ is a \Fr
space with group action induced by $\Kl$ on $\KsgX,$ $s\in\R.$
Globally on $\M$ we then define $\WW^{s,\gamma}_\P(\M)$ to be the
subspace of $\WW^{s,\gamma}(\M)$ locally near the edge described
by (\ref{WP}).\\

\nt Note that there is a slightly modified global edge calculus,
based on the spaces
$$\K^{s,\gamma;g}(X^\wedge):=\ang{r}^{-g}\KsgX$$

\nt for $g=s-\gamma$ rather than $\KsgX,$ and the group action
$$\kappa_\lambda^g:u(r,x)\mapsto\lambda^{g+\frac{n+1}{2}}u(\lambda
r,x),$$

\nt$\lambda\in\R_+,$ instead of (\ref{kl}). Edge spaces
$\WW^s(\R^q,E)$ modelled on $E=\K^{s,\gamma;g}(X^\wedge)$ with
alternative group actions of that kind have been suggested in
\cite{Aira2}, in connection with the solvability of hyperbolic
equations. In addition, following a remark of \cite{Tark4}, the
resulting global edge spaces on $\W$ have then particularly natural
invariance properties. All considerations here in connection with
§Green symbols and Green operators easily generalise to the modified
spaces. For simplicity we return to the case without $g;$ then in
order to have invariance of the global objects we need to choose an
atlas on $\W$ with some specified behaviour of transition maps.

\nt A Green operator $G$ (of the type of an upper left corner) with
(discrete or continuous) asymptotics is an operator that is locally
near $Y$ in stretched coordinates $(r,x,y)$ of the form
$$\Op_y(g)u(y)=\Int_{\R^q}\Int_{\R^q} e^{i(y-y')\eta}g(y,\eta)u(y')dy'\dslash\eta,$$
for a Green symbol $g(y,\eta)$ of order $\mu,$ modulo a global
smoothing operator. The latter category of operators is
characterised by the property to define a continuous map
$\WW^{s,\gamma}(\M)\rightarrow\WW^{\infty,\gamma-\mu}_\P(\M)$ for
some asymptotic type $\P$ and a similar property of the formal
adjoint. The global definition of Green operators is justified by
the following remark.
\begin{rem}
With symbols $g(y,\eta)\in S^\mu_{(\clas)}(\Omega\times\R^q;E,\til
E)$ in the sense of the notation in $\mathrm{Section}$ $\ref{1.2}$
we can form associated operators
$$\Op_y(g)u(y)=\Int_{\R^q}\Int_\Omega e^{i(y-y')\eta}g(y,\eta)u(y')dy'\dslash\eta,$$
$\dslash\eta:=(2\pi)^{-q}d\eta.$ In particular, if $g(y,\eta)$ is a
Green symbol of the kind $(\ref{gcont}),$ then for every
$\varphi(r)\in C_0^\infty(\R_+)$ the operators $\varphi\Op_y(g)$ and
$\Op_y(g)\varphi$ are smoothing on $\R_+\times\R^q\times X.$ In
particular, we see that the singularities of Green operators are
concentrated on the boundary $\set{0}\times\Omega\times X$ as is
expected in analogy to a corresponding behaviour of Green operators
in classical boundary value problems.
\end{rem}

\nt In fact, let us first note that the operators of multiplication
$\MM_\varphi$ by $\varphi$ generate (non-classical) symbols
$\MM_\varphi\in S^0\Set{\Omega\times\R^q;\KsgX,\KsgX}$ for every
$\gamma\in\R.$ Moreover, the multiplication by $r^{-N}\varphi$ for
any $N\in\N$ is of similar behaviour. Then,
$$\varphi g(y,\eta)=r^{-N}\varphi r^Ng(y,\eta)=r^{-N}\varphi[\eta]^{-N}(r[\eta])^Ng(y,\eta).$$
Since the order of $(r[\eta])^Ng(y,\eta)$ is the same as that of
$g(y,\eta)$ and the multiplication by $[\eta]^{-N}$ gives rise to an
order shift by $-N$ we obtain that $\varphi g(y,\eta)$ is an
operator-valued symbol of order $-\infty,$ and hence
$\varphi\Op_y(g)=\Op_y(\varphi g)$ is smoothing.\\

\nt Observe that a Green operator $G$ on $\M$ induces continuous
operators
$$\WW^{s,\gamma}(\M)\rightarrow\WW^{s-\mu,\gamma-\mu}_\P(\M)$$
for every $s\in\R,$ where $\P$ is a (discrete or continuous)
asymptotic type associated with $G.$ This is a consequence of
general continuity on Sobolev spaces.
\begin{rem}
Green operators on a $($stretched$)$ manifold $\M$ with edges form
an algebra, and the composition is compatible with the local
symbolic structure; in particular, the homogeneous principal symbols
$($in the sense of twisted homogeneity$)$ behave multiplicatively.
\end{rem}

\subsection{Green operators with parameters}

The concept of operator-valued symbols as in Section \ref{1.2} has a
parameter-dependent analogue, when we replace the covariable
$\eta\in\R^q$ by $(\eta,\lambda)\in\R^q\times\R^l$ and require the
symbolic estimates with respect to $(\eta,\lambda).$ In particular,
we obtain a generalisation of Definition \ref{Green symbol} to the
$\lambda$-dependent case, cf. also Remark
\ref{3.1.2}.\\

\nt The construction of the preceding section then gives us
parameter-dependent families of Green operators. According to the
iterative concept of building up pseudo-differential calculi on
manifolds with higher (polyhedral) singularities we may employ such
parameter-dependent families as (operator-valued) symbols of a next
generation of operators, for instance, on the infinite (stretched)
cone $\R_+\times\M$ with base $\M.$ Constructions in that sense may
be found in the paper \cite{DeDo1}, in particular,
a number of kernel cut-off results for such operator functions.\\

\nt Kernel cut-offs can be organised on the level of symbols. In
order to illustrate the effects we want to consider Green symbols
with discrete asymptotics as in Definition \ref{Green symbol} which
belong to spaces of the kind
\begin{equation}\label{parcl}
S^\mu_\clas\big(\Omega\times\R^q\times\R^l;\KsgX,\S^\beta_\P(X^\wedge)\big)
\end{equation}
such that $g^*(y,\eta,\lambda)$ belongs to corresponding analogue of
the space in (\ref{gstar}). For $l=1$ we also write $\Gamma_\delta$
instead of $\R$ when $\lambda$ is involved in the form
$z=\delta+i\lambda$ for some $\delta\in\R.$ Moreover, let
$$S^\mu_\clas\big(\Omega\times\R^q\times\C;\KsgX,\S^\beta_\P(X^\wedge)\big)
$$
denote the space of all $g(y,\eta,z)$ which are holomorphic in
$z\in\C$ such that $g(y,\eta,\delta+i\lambda)$ belongs to
(\ref{parcl}) (for $l=1)$ for every $\delta\in\R,$ uniformly in
compact $\delta$-intervals, and where $g^*(y,\eta,z)$ satisfies an
analogous condition.
\begin{thm}\label{cut}
For every $\delta\in\R$ there is a continuous map
$$S^\mu_\clas\big(\Omega\times\R^q\times\Gamma_\delta;\KsgX,\S^\beta_\P(X^\wedge)\big)
\rightarrow
S^\mu_\clas\big(\Omega\times\R^q\times\C;\KsgX,\S^\beta_\P(X^\wedge)\big),$$
$g(y,\eta,\delta+i\lambda)\mapsto h(y,\eta,\delta+i\lambda),$ such
that
$$g(y,\eta,z)-h(y,\eta,z)|_{\Omega\times\R^q\times\Gamma_\delta}\in
S^{-\infty}\big(\Omega\times\R^q\times\Gamma_\delta;\KsgX,\S^\beta_\P(X^\wedge)\big)$$
for every $s\in\R$ $($and such that the pointwise formal adjoint
have an analogous property$).$
\end{thm}

\nt Theorem \ref{cut} can be proved by applying a kernel cut-off
argument as used in an analogous context in \cite{DeDo1}. A similar
result holds for Green symbols with continuous
asymptotics.\\

\nt According to Theorem \ref{THM2.1.1} the holomorphic symbol
$h(y,\eta,z)$ has a family of integral kernels
\begin{equation}\label{f}
f_{\Re z}(r[\eta,\lambda],x,r'[\eta,\lambda],x',y,\eta,\lambda),
\end{equation}
$\lambda=\Im z,$ via a representation of the form (\ref{gg}). This
is valid for every fixed $\Re z;$ however, the holomorphic
dependence of (\ref{f}) on $z$ is by no means obvious. In other
words, kernel cut-off constructions which produce a holomorphy in a
complex covariable are better applied to the symbols in their
original definition rather than their integral kernels.\\


\bibliographystyle{wileybib}
\bibliography{master}

\end{document}